\documentclass[a4paper,10pt]{article}
\usepackage{amsfonts, amsthm, amsmath, amssymb, graphicx, color}
\usepackage{tikz}
\usepackage[normalem]{ulem}

\title{On continued fraction maps associated with a submodule of $\mathfrak o(\sqrt{-3})$} 
\author{
Hitoshi Nakada \\
Department of Mathematics, Keio University, \\ 
Yokohama, Japan \\
nakada@math.keio.ac.jp \\[5pt]
Rie Natsui \\
Department of Mathematics, Physics and Computer Science,\\
 Japan Women's University\\
Tokyo, Japan\\
natsui@fc.jwu.ac.jp \\[5pt]
Mako Toyosumi \\
Department of Mathematics, Physics and Computer Science, \\
Japan Women's University\\
Tokyo, Japan \\
m1716066tm@ug.jwu.ac.jp
}

\newtheorem{thm}{Theorem}

\newtheorem{lem}{Lemma}
\newtheorem{prop}{Proposition}
\newtheorem{cor}{Corollary}
\newtheorem{rem}{Remark}
\newtheorem{defn}{Definition}

\newcommand{\confrac}[2]{
\frac{\displaystyle{
\strut\hfill{#1}\hfill\;\vrule}}
{\displaystyle{
 \strut\vrule\;\hfill{#2}\hfill}}}

\begin{document}
\maketitle
\begin{abstract}
We define a continued fraction map $T$ associated with the 
$\mathfrak o(\sqrt{-3})$-module $\mathcal J = \eta \cdot\mathfrak o(\sqrt{-3})$, 
$\eta = \frac{3 + \sqrt{-3}}{2}$, which is an Eisenstein field version of the continued 
fraction map associated with  $\mathfrak o(\sqrt{-1}) \cdot (1 + i)$ defined by J.~Hurwitz in the case of the Gaussian field.  Together with $T$, we show that all complex numbers $z$ can be expanded as $\mathcal J$-coefficients.  
We discuss some basic properties of these continued fraction expansions such as 
the monotonicity of the absolute value of the denominator of the principal convergent $q_{n}$ and the existence of the absolutely continuous ergodic invariant probability measure for $T$.    
\\[3pt]
\\[7pt]
\end{abstract}
\section{Introduction}
Let $d \ge 1$ be a square free positive integer and $\mathfrak o(\sqrt{-d})$ be the 
set of algebraic integers of $\mathbb Q(\sqrt{-d})$, i.e., 
\[
\mathfrak o(\sqrt{-d}) = \left\{\begin{array}{lcc}
\left\{ n\,\frac{1 + \sqrt{-d}}{2} + m\,\frac{1 - \sqrt{-d}}{2} : n, m \in 
\mathbb Z\right\} & \mbox{if} & d \equiv 3 \,\,\mbox{mod} \,\, 4\\
\left\{n + m \sqrt{-d} : n, m \in \mathbb Z \right\} & \mbox{otherwise} &{}
\end{array}\right.
\]
Here we only consider the case $d = 1, 2, 3, 7, 11$, which means that 
the Euclidean algorithm works for $\mathfrak o(\sqrt{-d})$.  For those $d$, there 
exist the nearest integer type complex continued fraction maps with the 
$\mathfrak o(\sqrt{-d})$-coefficients for each $d$, see \cite{La},  
\cite{E-I-N-N}, \cite{E-N-N-1}, and \cite{E-N-N-2}.  
On the other hand, J.~Hurwitz \cite{J-Hu} defined 
a nearest integer type complex continued fraction map with the 
$(1 + i)\cdot \mathfrak o(\sqrt{-1})$-coefficients.  Then S.~Tanaka \cite{Ta}, N.~Oswald \cite{Os} and H.~Nakada \cite{N} discussed its ergodic properties.  
J.~Hurwitz's map can be regarded as a complex variant of the continued fractions 
with even partial coefficients for real numbers, see Schweiger \cite[p.18]{Sch} for the real case.  There is a significant difference between 
the real case and the complex case in ergodic point of view.  In the case of real 
numbers, the associated map has an absolutely continuous invariant measure but 
it is not a finite measure, on the other hand,  in the complex case, there exists an absolutely continuous invariant probability measure.  The difference implies that 
$\lim_{n \to \infty}\tfrac{1}{n} \log |q_{n}|$ has a positive finite limit (a.e.) 
in the case of the complex numbers, see \S4, on the other hand, the limit is $0$ 
(a.e.) the real case case, where $q_{n}$ denotes the denominator of the 
$n$th convergent of a continued fraction.        

In this point of view, it would be natural to consider this type of complex continued 
fractions for other imaginary quadratic fields.  It is possible to define the nearest integer type map associated with a submodule of $\mathfrak o(\sqrt{-d})$,
$d=1, 2, 3, 7, 11$.  However, except for the case $d=3$ with the submodule 
$\eta\cdot \mathfrak o(\sqrt{-3}) = \{\eta \cdot a : a \in \mathfrak o(\sqrt{-3})\}$, 
$\eta = \frac{3 + \sqrt{-3}}{2}$, for each map, 
there always exists a set $A$ of positive (complex) Lebesgue measure such that the 
complex continued expansion of $z \in A$ never converges to $z$ itself.   
This fact follows the same as the nearest type continued fraction maps with 
non-Euclidean imaginary quadratic case, see \cite[\S6]{E-N-N-1}.  

The aim of this paper is to give some fundamental properties of the continued 
fractions associated with $\mathcal J = \eta\cdot\mathfrak o(\sqrt{-3}) $, 
$\eta = \frac{3 + \sqrt{-3}}{2}$.  First, we will see that for any $z \in U$,  
there exists a continued fraction expansion of $z$ with $\eta \cdot
\mathcal J$-quotients, where $U$ is the fundamental domain of 
$\mathbb C / \mathcal J$ defined below. 
This expansion is given by the following map $T$ except for 
$z = \overline{\zeta}$ and $-\zeta$.  We have to give expansions for 
these two points separately.   

\begin{figure} \label{fig-1-1}
\begin{center}
\begin{tikzpicture}[scale=0.4]
\fill[lightgray](4,0)--(2, 3.464)--(-2, 3.464)--(-4,0)--(-2, -3.464)--(2, -3.464)--(4, 0)--
cycle;
\draw[dashed] (4,0)--(2, 3.464)--(-2, 3.464)--(-4,0)--(-2, -3.464)--(2, -3.464)--(4, 0)--
cycle;

\draw[thick] (2, -3.464)--  (4,0); 
\draw[thick] (-2, 3.464)-- (2, 3.464); 
\draw[thick] (-2, -3.464)--  (-4,0);

\draw (2, 3.464)-- (0, 0) --(2, -3.464); 
\draw (0, 0) --(-4, 0); 

\draw (0,0) circle(4);

\node at (0, 0) {\Large $\cdot$}; 
\node at (-2, -3.464) {\Huge$\cdot$}; 
\node at (2, -3.464) {\Huge $\cdot$}; 

\node at (0, -0.7) {\large $0$}; 
\node at (4.3, -0.3) {$1$}; 
\node at (-4.5, -0.3){$-1$};
\node at (3.9, 4) {$\zeta = \frac{1 + \sqrt{-3}}{2}$} ;
\node at (-2.8, -4.0) {$-\zeta$} ;
\node at (2.2, -4.0) {$\overline{\zeta}$}; 

\node at (2, 0) {$\underbar{U}$}; 
\node at (-1, 1.7) {$\zeta^{2}\underbar{U}$}; 
\node at (-1, -1.7) {$\zeta^{-2}\underbar{U}$}; 

\node at (0, -6) {$U$: grayed area};

\end{tikzpicture}
\qquad \qquad
\begin{tikzpicture}[scale = 0.4]
\fill[lightgray] (4, 0) arc (-30:90:2.31) arc (30:150:2.31) arc (90:210:2.31) arc 
(150:270:2.31) arc (210:330:2.31) arc (-90:30:2.31);
\draw (4,0)--(2, 3.464)--(-2, 3.464)--(-4,0)--(-2, -3.464)--(2, -3.464)--(4, 0)-- cycle;
\draw (4, 0) arc (-30:90:2.31) arc (30:150:2.31) arc (90:210:2.31) arc (150:270:2.31) 
arc (210:330:2.31) arc (-90:30:2.31);

\node at (0,-0.4) {\large $0$}; 
\node at (4.3, -0.3) {$1$}; 

\node at (0, -6) {$U^{-1}$: outside of grayed area};

\end{tikzpicture}

\end{center}
\caption{$U$ and $U^{-1}= \{z : \tfrac{1}{z} \in U\}$}
\end{figure}
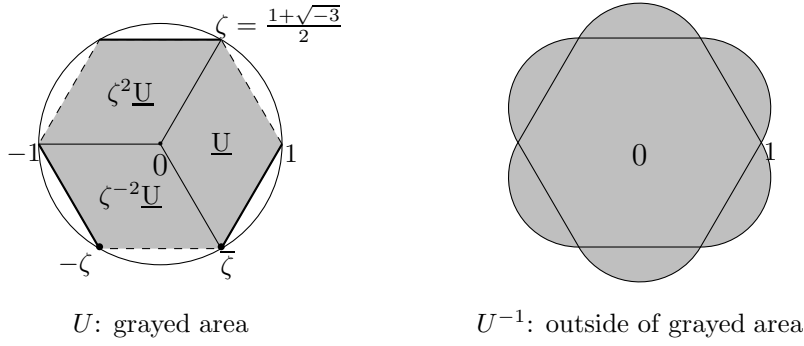

We put 
\[
U := \underbar{U} \cup \zeta^{2}\underbar{U} \cup \zeta^{-2}\underbar{U} 
\cup  \{-\zeta, \overline{\zeta}\} 
\]
with $\underbar{U} := \{ z = x \zeta + y \overline{\zeta} : 0\le x<1,\, 
0\le y \le$ 1 \}, $\zeta = \frac{1 + \sqrt{-3}}{2}$,  see Fig.\ref{fig-1-1}. 
Here and henceforth, 
for $z \in \mathbb C$ and $A \subset \mathbb C$, $\overline{z}$, $A+ z$, $z A$ 
denote the complex conjugate of $z$, $\{w + z : w \in A\}$, $\{zw : w \in A\}$, 
respectively.  
\begin{rem} There are other choices for three sides and two corners of the hexagon.   
However, there is no serious difference among choices. 
\end{rem}
  
It is easy to see that $U$ is the fundamental region of $\mathbb C$:
\[
\bigcup_{a \in \mathfrak o(\sqrt{-3})} a + U = \mathbb C \quad \mbox(disj) . 
\]
Thus, for every $z \in \mathbb C$, there exists a unique 
$a \in \mathfrak o(\sqrt{-3})$ such that $z \in a + U$.  We denote by $[z]$ for this 
$a$.   We define the continued fraction map $T$ by 
\begin{equation} \label{eq-1-1}
T(z) = \left\{ \begin{array}{ccc} 
 \frac{1}{z} - \left[ \frac{1}{z} \right] & \mbox{if} & z \in U \setminus \{0\} \\
0 & \mbox{if} & z = 0
\end{array} \right.
\end{equation}
As usual, we put $b_{n}(z) = \left[ \frac{1}{T^{n-1}(z)}\right]$ with 
$\tfrac{1}{0} = \infty$, $\tfrac{1}{\infty} = 0$ and and $[\infty] = \infty$ 
and get 
\begin{equation} \label{eq-1-2}
z = \confrac{1}{b_{1}(z)} + \confrac{1}{b_{2}(z)} + \cdots +
\confrac{1}{b_{n}(z) + T^{n}(z)}
\end{equation} 
with $\tfrac{1}{\infty} = 0$. 
In \S2, we discuss the convergence of the continued fraction expansion, namely 
\[
 \lim_{n \to \infty} \confrac{1}{b_{1}(z)} + \confrac{1}{b_{2}(z)} + \cdots + 
\confrac{1}{b_{n}(z)} .
\]
Then we discuss some properties of continued fractions associated with $T$.  
We denote by 
$\tfrac{p_{n}(z)}{q_{n}(z)}$ the $n$th convergent of $z$, see \S2 for its explicit definition.  
 
In \S3, we show that $T$ has the finite range structure, see \cite{Y}.  
This property shows the existence of the absolutely continuous invariant measure for $T$.  Moreover, in \S4, we construct a sort of ``dual area" of $U$, which induces a representation of the natural extension of $T$.  
It would be interesting that the boundary of this ``dual area" is a union of arcs.  
This is somehow similar to the J. Hurwitz case, see \cite{N, Ta}.    
This fact is remarkable if one compares to the nearest integer type continued fraction 
case, where the boundaries consist of ``fractal curves", see \cite{E-I-N-N, E-N-N-2}.  

We also see that the absolutely continuous invariant measure for $T$ is finite.  
From this construction, we can show that $|q_{n}(z)|$, $ n\ge 1$, is strictly 
monotone increasing for any $z \in U$.   For almost every $z \in U$, it is a direct 
consequence of \S4.  
However, to show the monotonicity for any $z \in U$, we need 
careful discussion for $z$ such that $T^{n}(z) = -\zeta$ or 
$\overline{\zeta}$, $n\ge 0$ and we do it in \S5.   

%%%%%%%%%%%%%%%%%%%%%%%%%%%%
%%%%%%%%%%%%%%%%%%%%%%%%%%%%
\section{Convergence of continued fractions}
In the sequel, we always assume that $\eta = \frac{3 + \sqrt{-3}}{2}$.  
As usual, we put 
\begin{equation} \label{eq-2-1}
\begin{pmatrix} p_{n-1}(z)& p_{n}(z) \\ q_{n-1}(z) & q_{n}(z) \end{pmatrix} := 
\begin{pmatrix} 0& 1 \\ 1 & b_{1}(z) \end{pmatrix}
\begin{pmatrix} 0& 1 \\ 1 & b_{2}(z) \end{pmatrix}
\cdots 
\begin{pmatrix} 0& 1 \\ 1 & b_{n}(z) \end{pmatrix}
\end{equation} 
for $n \ge1$, which implies 
\begin{equation} \label{eq-2-2}
\frac{p_{n}(z)}{q_{n}(z)} = 
\confrac{1}{b_{1}(z)} + \confrac{1}{b_{2}(z)} + 
\cdots + \confrac{1}{b_{n}(z)} , \,\,\,  p_{n}(z)\,\, \mbox{and}\,\, q_{n}(z)\,\,
\mbox{are coprime}. 
\end{equation}  
For a given sequence $(b_{n} : b_{n} \in \mathcal J , \, n\ge 1)$, we also write 
\[
\begin{pmatrix} p_{n-1}& p_{n} \\ q_{n-1} & q_{n} \end{pmatrix}  = 
\begin{pmatrix} 0& 1 \\ 1 & b_{1} \end{pmatrix}
\begin{pmatrix} 0& 1 \\ 1 & b_{2} \end{pmatrix}
\cdots 
\begin{pmatrix} 0& 1 \\ 1 & b_{n} \end{pmatrix} .  
\] 
We will show that for any $z \in U$,
\begin{equation}  \label{eq-2-3}
\lim_{n \to \infty} \frac{p_{n}(z)}{q_{n}(z)} = z  
\end{equation}
if $T^{n}(z) \ne \overline{\zeta}$ nor $-\zeta$ for $n\ge 1$.  
For each $2 \times 2$ matrix, we consider the linear fractional map defined by the matrix.  We denote it also by the same matrix, i.e. 
$\begin{pmatrix}a & b \\ c & d\end{pmatrix}(z) = \tfrac{az+ b}{cz + d}$.  
 
We put $z_{n}= T^{n}(z)$, $n \ge 0$ with $z_{0} = T^{0}(z) = z$, and get 
\[
z =  \begin{pmatrix} p_{n-1}(z)& p_{n}(z) \\ q_{n-1}(z) & q_{n}(z) \end{pmatrix}
( z_{n} )
= 
\frac{p_{n-1}(z) z_{n} +  p_{n}(z)}{q_{n-1}(z) z_{n} +  q_{n}(z)}.  
\]
There are three cases: (i) $z_{n} \notin \{0, \overline{\zeta}, -\zeta\}$ for any 
$n \ge 0$, 
(ii) $z_{n} = 0$ for some $n \ge 0$,
(iii) $z_{n} =\overline{\zeta}$ or $ -\zeta$ for some $n \ge 0$.  

In the case (ii), it is easy to see that 
\[
\begin{pmatrix} p_{m-1}(z)& p_{m}(z) \\ q_{m-1}(z) & q_{m}(z) \end{pmatrix}
= 
\begin{pmatrix} p_{n-1}(z)& p_{n}(z) \\ q_{n-1}(z) & q_{n}(z) \end{pmatrix}
\]
and $z_{m}=0$ for $m \ge n$.  Thus we have $z = \tfrac{p_{m}(z)}{q_{m}(z)}$ 
for $m \ge n$ and thus \eqref{eq-2-3} holds. 
To consider the case (iii), we have to redefine $b_{n+k}(z)$, $k \ge 0$, which we will consider later. 

To prove \eqref{eq-2-3} for cases (i) and (iii), we use \cite[Prop. 2-1]{D}.  Indeed, 
we use Prop. 2-1 (i) and (iii) for our cases (iii) and (i) respectively.  We use  
\begin{equation}  \label{eq-2-4}
 \left|z - \frac{p_{n}(z)}{q_{n}(z)}\right|  = \left|\frac{1}{q_{n}(z)}\right| \cdot 
|q_{n}(z)\cdot z - p_{n}(z)| = \left|\frac{1}{q_{n}(z)}\right| \cdot | z_{0} \cdot z_{1} 
\cdots z_{n}|, \,\, n \ge 1 , 
\end{equation}
which follows by induction. With this estimate, we only need to show either 
$\lim_{n \to \infty} |z_{0}z_{1} \cdots z_{n}| = 0$ or $\lim_{n\to\infty}|q_{n}(z)| = \infty$, since $q_{n}(z) \in \mathfrak o(\sqrt{-3})$ and $|z_{m}| \le 1$ 
for $m \ge 0$.  We note that $q_{n}(z) \ne 0$, $n\ge0$,  in our cases (i) and (iii).   

If $\left|q_{n}(z)\right|$ is strictly monotone increasing, then it is easy to see 
\eqref{eq-2-3}.  Indeed, the following well-known formula 
\begin{align} \label{eq-2-5} 
\left|z - \frac{p_{n}(z)}{q_{n}(z)}\right| &
= 
\left| \frac{p_{n-1}(z) z_{n} +  p_{n}(z)}{q_{n-1}(z) z_{n} +  q_{n}(z)} - 
\frac{p_{n}}{q_{n}}\right| \notag \\
{} &\le  \left|\frac{z_{n}}{q_{n}(z) (q_{n-1}(z) z_{n} + q_{n}(z))}\right|  \notag\\
{} & \le 
\frac{1}{\left| q_{n}(z) (q_{n-1}(z) z_{n} + q_{n}(z)) \right| }.     
\end{align}
leads us the conclusion, since $\left|(q_{n-1}(z) z_{n} + q_n(z)) \right| > \sqrt{3}-1$  
for $z \in U$ in this case.  However, the proof of the monotonicity is not simple 
and we postpone it to the following sections. 
Anyway we have the following with \eqref{eq-2-4} for the case (i).  For the case (iii), 
we give $(b_{n})$ in a different way. 
%%%%%%%%%%%%%%%%%%%%%
\begin{thm} \label{thm-2-1}
For any $z \in U$, there exists a sequence $(b_{n}:n \ge 1)$ such that $b_{n} \in 
\mathcal J$, $n \ge 1$, and 
\[
z = \lim_{n \to \infty} \frac{p_{n}}{q_{n}} 
\]
\[
\left\{ \begin{array}{lll} 
p_{n+1} = b_{n+1}p_{n} + p_{n-1} \\
q_{n+1} = b_{n+1}q_{n} + q_{n-1} 
\end{array}\right.  
\]
for $n\ge 0$ and 
\[
\left\{ \begin{array}{lll} 
p_{-1}=1 , p_{0}=0, \\
q_{-1}=0 , q_{0}=1
\end{array}\right.  . 
\]
\end{thm}
%%%%%%%%%%%%%%%%%%%%%%%%%%%%%%

In the cases (i) and (ii), we can use $b_{n} = b_{n}(z)$ by $T$.   
The only difficulty of the proof of this theorem is the case (iii).   
Indeed, if we follow the definition of $b_{n}(z)$, then $b_{n}(-\zeta) = \sqrt{-3}$ for any $n \ge 1$.  It is easy to see that 
\[
\underbrace{\confrac{1}{\sqrt{-3}} + \confrac{1}{\sqrt{-3}} + \cdots + \confrac{1}{\sqrt{-3}}}_{6} = 0, 
\]
which shows that $\lim_{n \to \infty}\confrac{1}{\sqrt{-3}} + 
\confrac{1}{\sqrt{-3}} + \cdots$ does not converge.  
For this reason,  we have to re-define $b_{n+k}(z)$, $k \ge 1$ 
if $T^{m}(z) \ne \{\overline{\zeta}, -\zeta\}$ for $0 \le m < n$ and 
$T^{n}(z) \in \{\overline{\zeta}, -\zeta\}$, $n\ge 0$.  
%%%%%%%%%%%%%%%
\begin{lem} \label{lem-2-1}
We have the expansion of $-\zeta$ and $\overline{\zeta}$ by 
\begin{equation} \label{eq-2-6}
{\rm (i)} \,\,\, -\zeta = \dot{\confrac{1}{\sqrt{-3}}} + \confrac{1}{\sqrt{-3}} +   
\confrac{1}{-\overline{\eta}} + \dot{\confrac{1}{\eta}} +  \cdots , 
\end{equation}
and 
\begin{equation} \label{eq-2-7}
{\rm (ii)} \,\,\overline{\zeta}  = \dot{\confrac{1}{\sqrt{-3}}} + \confrac{1}{\sqrt{-3}} 
+ \confrac{1}{\eta} + \dot{\confrac{1}{-\overline{\eta}}} +  \cdots , 
\end{equation}
respectively, 
i.e.,  
\[
b_{4n+1}(-\zeta)  = b_{4n+2}(-\zeta) = \sqrt{-3},  \,\,
b_{4n+3}(-\zeta)= - \overline{\eta}, \,\, b_{4n+4}(-\zeta) = \eta 
\]
and
\[
b_{4n+1}(\overline{\zeta})  = b_{4n+2}(\overline{\zeta}) = \sqrt{-3}, \,\,\,  
b_{4n+3}(\overline{\zeta})= - \eta, \,\,\,b_{4n+4}(\overline{\zeta}) = \overline{\eta}
\]
for $n \ge 0$.
\end{lem}
\begin{rem} The explanation of finding (i) and (ii) will be discussed in the proof of 
Theorem \ref{thm-5-1}.   
\end{rem} 
{\bf Proof.} 
We put 
\begin{equation} \label{eq-add}
\begin{array}{lll}
S_{-\zeta, 0} & = &\left\{z : {\rm Im}\, z = - \frac{\sqrt{3}}{2}, \,\, {\rm Re} \, z 
< -\tfrac{1}{2}\right\}\cup \{\infty\}\\
S_{-\zeta, 1} & = &\left\{z : \left| z + \frac{2\sqrt{-3}}{3}\right| = 
\frac{\sqrt{3}}{3}, \, {\rm Im} \, z < -\tfrac{\sqrt{3}}{2}, \, {\rm Re} \, z \le 0 
\right\} \\
S_{- \zeta, 2} & = &\left\{z : \left|z + \left( \tfrac{\sqrt{-3}}{3}\right) \right| 
= \frac{\sqrt{3}}{3}, \, 
{\rm Im} z < \frac{-\sqrt{3}}{2}, \, {\rm Re} \, z \le 0 \right\} \cap \{z:|z| > 1\}\\
S_{-\zeta, 3} & = &\left\{z :  {\rm Im} z = 0, \,\, {\rm Re}\, z >1 \right\}
\end{array} . 
\end{equation}
and 
\[
\begin{array}{lll}
S_{\overline{\zeta}, 0} & = &\left\{z : {\rm Im}\, z = - \frac{\sqrt{3}}{2}, \,\, {\rm Re} \, z 
> \frac{1}{2}\right\}\cup \{\infty\}\\
S_{\overline{\zeta}, 1} & = &\left\{z : \left| z + \frac{2\sqrt{-3}}{3}\right| = \frac{\sqrt{3}}{3}, 
\, {\rm Im} \, z < \tfrac{\sqrt{3}}{2}, \, {\rm Re} \, z \ge 0 \right\} \\
S_{\overline{\zeta}, 2} & = &\left\{z : \left|z - \left(- \tfrac{\sqrt{-3}}{3}\right) \right| = \frac{\sqrt{3}}{3}, \, 
{\rm Im} z < \frac{-\sqrt{-3}}{2}, \, {\rm Re} \, z \ge 0 \right\} \\
S_{\overline{\zeta}, 3} & = &\left\{z :  {\rm Im} z = 0, \,\, {\rm Re}\, z < -1 \right\}
\end{array} . 
\]
We start with $(z_{0}, S_{0}) = (-\zeta, S_{-\zeta, 0})$ and apply the map 
$\tfrac{1}{z} - \sqrt{-3}$, $z \in \mathbb C$, for both components.  Then we have 
$(z_{1}, S_{1}) = (-\zeta, S_{-\zeta, 1})$.  Again we apply the same map and have 
$(z_{2}, S_{2}) = (-\zeta, S_{-\zeta, 2})$.  Next we apply $\tfrac{1}{z} - 
\left(-\overline{\eta}\right)$ and have 
$(z_{3}, S_{3})$ which is a subset of $(1, S_{-\zeta, 3})$.  We continue by $\tfrac{1}{z} - 
\eta$ and have $(z_{4}, S_{4}) \subset 
(-\zeta, S_{-\zeta, 0})$.  We continue this process inductively and get a sequence 
\[
\left (z_{4n + k}, S_{4n + k}) \,\, \mbox{such that} \,\, (z_{4n + k}, S_{4n+k}) 
\subset (z_{k}, S_{-\zeta, k}) \,\,\mbox{for} \,\, n \ge 0, 0 \le k \le 3\right. .
\] 
We note that for any $\omega \in S_{-\zeta, k}, 0 \le k \le 3, |\omega| >1$.Now we consider $\tfrac{q_{0}}{q_{-1}} = 1/0 = \infty$.  Then $\left( \tfrac{q_{0}}{q_{-1}}
- \sqrt{-3}\right)^{-1} = \tfrac{q_{1}}{q_{0}} = -\sqrt{-3} \in S_{1}$.  Inductively, we 
see $\tfrac{q_{k}}{q_{k-1}} \in S_{k}$, which shows 
$|q_{k}| > |q_{k-1}|$ for $k \ge 1$.  In this case, we have $|z_{k}| = 1$ for 
$k \ge 0$.  According to \eqref{eq-2-4}, we get \eqref{eq-2-3} for $-\zeta$.  
It is easy to see that the convergence also follows for $\overline{\zeta}$ with 
$S_{\overline{\zeta},k}$, $0 \le k \le 3$ by the same way.  
\qed 
%%%%%%%%%%%%%%%%%%%%

{\bf Proof of Theorem \ref{thm-2-1}.} 
We already proved the case $T^{n}(z) = 0$ for some $n \ge 0$. 
From Lemma \ref{lem-2-1} and \eqref{eq-1-2}, 
we have the expansion of $z \in U$ such that $T^{n}(z) = -\zeta$ or 
$\overline{\zeta}$ for some $n\ge 0$ as   
\[
z = \confrac{1}{b_{1}(z)} + \confrac{1}{b_{2}(z)} + \cdots + \confrac{1}{b_{n}(z)} + 
\confrac{1}{b_{1}(-\zeta)} + \confrac{1}{b_{2}(-\zeta)} + \cdots 
\]
or 
\[
z = \confrac{1}{b_{1}(z)} + \confrac{1}{b_{2}(z)} + \cdots + \confrac{1}{b_{n}(z)} + 
\confrac{1}{b_{1}(\overline{\zeta})} + \confrac{1}{b_{2}(\overline{\zeta})} + \cdots , 
\]
respectively. 
If $T^{n}(z) \ne 0$, $-\zeta$ nor $\overline{\zeta}$ for any $n \ge 0$, 
then  $|z| < 1$ for any $z \in U\setminus \{- \zeta, \overline{\zeta}\}$ and $|T'(z)|> 1$, we note that $\zeta, -\overline{\zeta}, \pm1\notin U$. 
So, for sufficiently small $\varepsilon >0$, $|z_{n}|< 1 -\varepsilon$ holds 
for infinitely many $ n \ge0$.  Hence, $|z_{0} z_{1} z_{2} \cdots z_{n}| \to 0$ as 
$n \to \infty$.  From \eqref{eq-2-4}, we get the assertion of the theorem.  
\qed \\

%%%%%%%%%%%%%%%%%%%%
%%%%%%%%%%%%%%%%%%%%
%%%%%%%%%%%%%%%%%%%%
\section{Finite range structure}
For a sequence $(a_{1}, a_{2}, \ldots, a_{n}) \in \mathcal J^{n}$, we define the cylinder set of rank $n$ associated with 
$(a_{1}, a_{2}, \ldots, a_{n})$ by 
\[
\langle a_{1}, a_{2}, \ldots, a_{n} \rangle 
= \{z \in U : b_{1}(z) = a_{1}, b_{2}(z) = a_{2}, \ldots , b_{n}(z) = a_{n} \} .  
\]
If a cylinder set $\langle a_{1}, a_{2}, \ldots, a_{n} \rangle $ has an inner point, 
then we call $(a_{1}, a_{2}, \ldots, a_{n})$ an admissible sequence.  
If 
\[
\{\left(T^{n}(\langle a_{1}, a_{2}, \ldots, a_{n} \rangle)\right)^{\circ} :  
(a_{1}, \ldots, a_{n})\,\,\mbox{is admissible},\,\,n \ge 1\} 
\]
is a finite set, then $T$ is said to have a finite range structure, where 
$A^{\circ}$ denotes the interior of $A \subset \mathbb C$.  We will show that 
$T$ has this property. 

We put $U_{0} = U^{\circ}$ and 
\[
\begin{array}{lll}
U_{1,k} &= &\zeta^{k-1}\left(U_{0}\cap 
\left\{z : \left|z - \left(-\tfrac{2}{3} \eta\right)\right|> \tfrac{\sqrt{3}}{3}
\right\}\right) \\
U_{2,k} &= &\zeta^{k-1}\left(U_{0}\cap 
\left\{z : \left|z - \left(-\tfrac{1}{3} \eta\right)\right|> \tfrac{\sqrt{3}}{3}
\right\}\right) \\
U_{3,k} &= &\zeta^{k-1}\left(U_{0}\cap 
\left\{z = x+yi: y > \sqrt{3} x \right\} \right) \\
U_{4,k} &= &\zeta^{k-1}\left(U_{1, k}\cap 
\left\{z = x+yi : y > \sqrt{3}x  \right\} \right) \\
U_{5,k} &= &\zeta^{k-1}\left(U_{1, k+1}\cap 
\left\{z = x+yi:  y < \sqrt{3} x\right\} \right) \,: \, k+1 \,\,
\mbox{means} \,\, k+1\, \mbox{(mod. 6)}
\end{array}
\] 
for $1 \le k \le 6$, see Fig.2.   
If $\left( T^{n}(\langle a_{1}, a_{2}, \ldots, a_{n} \rangle)\right)^{\circ} = U_{0}$, then $\langle a_{1}, a_{2}, \ldots, a_{n} \rangle$ is said to be full. 
We also put 
\[
\begin{array}{lll}
V_{1,\ell}& = & \zeta^{\ell -1}\left\{z \in U : \,\,\left|z - \tfrac{1}{3}\overline{\eta}
\right|< \tfrac{\sqrt{3}}{3}, \,\, 
\left|z - \tfrac{1}{3}\sqrt{-3}\right|< \tfrac{\sqrt{3}}{3}\right\} \\
%%%%%%
V_{2,\ell}& = & \zeta^{\ell -1}\left\{z = x + i y \in U : \,\,\left|z - \tfrac{2}{3}\eta \right|> \tfrac{\sqrt{3}}{3}, \,\, 
\left|z - \tfrac{1}{3}\sqrt{-3}\right|> \tfrac{\sqrt{3}}{3}, \,\, 
x > 0,\,\, y >0 \right\} \\
V_{3,\ell}& = & \zeta^{\ell -1}\left\{z = x + i y \in U : \,\,\left|z - 
\tfrac{2}{3}{\eta}\right|> \tfrac{\sqrt{3}}{3}, \,\, 
\left|z - \tfrac{1}{3}\overline{\eta}\right|> \tfrac{\sqrt{3}}{3}, \,\, 
y >0,\,\, y < \sqrt{3} x \right\}\\
V_{4,\ell}& = & \zeta^{\ell -1} \left\{z \in U : \,\,\left|z - \tfrac{1}{3}\overline{\eta}
\right|< \tfrac{\sqrt{3}}{3}, \,\, 
\left|z - \tfrac{2}{3}\eta\right|< \tfrac{\sqrt{3}}{3}\right\} \\
V_{5,\ell}& = & \zeta^{\ell -1}\left\{z \in U : \,\,\left|z - \tfrac{1}{3}\eta\right| < \tfrac{\sqrt{3}}{3}, \,\, 
\left|z - \tfrac{1}{3}\sqrt{-3}\right| < \tfrac{\sqrt{3}}{3}\right\} \\
	V_{6,\ell}& = & \zeta^{\ell -1}\left\{z \in U : \,\,\left|z - \tfrac{1}{3}\overline{\eta}\right|> \tfrac{\sqrt{3}}{3}, \,\, 
\left|z - \tfrac{1}{3}\sqrt{-3}\right|> \tfrac{\sqrt{3}}{3}, \,\,  x>0, 
\,\, y>0 \right\} 
\end{array}
\]
 for $1 \le \ell \le 6$, see Fig.\ref{fig-3-2}.  
We denote $\mathcal V = \{V_{k, \ell}: 1 \le k \le 6, 1 \le \ell \le 6\}$ the partition defined by $V_{k, \ell}$.  Equivalently, 
$\mathcal V$ is the partition induced from $U_{k, \ell}$, $1 \le k \le 5$, $1 \le \ell \le 6$. 

Since we discuss the dynamic behavior 
of $T$ in \S4,  we do not mention about the boundary of each $V_{k, \ell}$ for the 
moment.   However, we do need the behavior of $ z \in \partial V_{k, \ell}$ in 
\S5.  In the rest of this section, we 
only consider cylinder sets associated with admissible sequences.  So we only say 
``a cylinder set" without saying ``associated with an admissible sequence". 
%%%%%%%%%%%%%%%%%%%%%%
%%%%%%%%%%%%%%%%%%%%%%
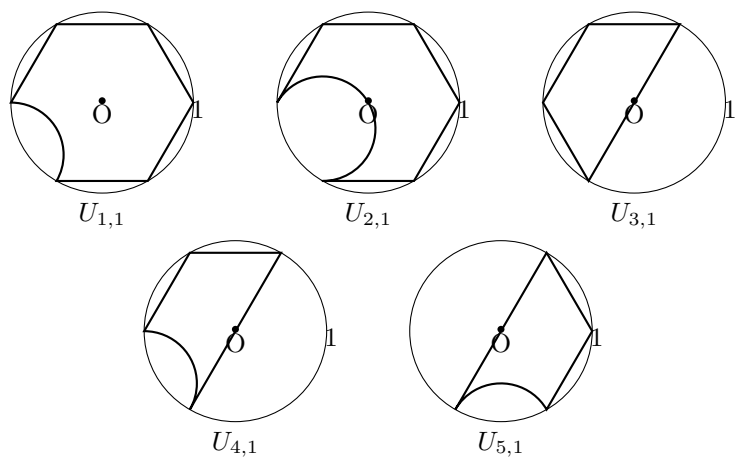
\begin{figure}      \label{fig-3-1}
\begin{center}
\begin{tikzpicture}[scale=0.3]
\draw[thick] (4,0)--(2, 3.464)--(-2, 3.464)--(-4,0) arc (90:-30:2.31) --
(2, -3.464)--(4, 0)--cycle;

\draw (0,0) circle(4);
\node at (0, 0) {\Huge $\cdot$}; 
\node at (0, -0.5) {${\rm O}$}; 
\node at (4.2, -0.3) {$1$}; 
\node at (0, -5) {$U_{1,1}$}; 
\end{tikzpicture}
\qquad
\begin{tikzpicture}[scale=0.3]
\draw[thick] (4,0)--(2, 3.464)--(-2, 3.464)--(-4,0) arc (150:-90:2.31)-- (-2, -3.464) -- 
(2, -3.464) -- (4, 0)--cycle;

\draw (0,0) circle(4);
\node at (0, 0) {\Huge $\cdot$}; 
\node at (0, -0.5) {${\rm O}$}; 
\node at (4.2, -0.3) {$1$}; 
\node at (0, -5) {$U_{2,1}$};

\end{tikzpicture}
\qquad 
\begin{tikzpicture}[scale=0.3]
\draw[thick] (2, 3.464)--(-2, 3.464)--(-4, 0)--(-2, -3.464) -- (2, 3.464) --cycle;

\draw (0,0) circle(4);
\node at (0, 0) {\Huge $\cdot$}; 
\node at (0, -0.5) {${\rm O}$}; 
\node at (4.2, -0.3) {$1$}; 
\node at (0, -5) {$U_{3,1}$};
\end{tikzpicture}

\begin{tikzpicture}[scale=0.3]
\draw[thick] (2, 3.464)--(-2, 3.464)--(-4, 0)arc (90:-30:2.31) -- (2, 3.464) --cycle;

\draw (0,0) circle(4);
\node at (0, 0) {\Huge $\cdot$}; 
\node at (0, -0.5) {${\rm O}$}; 
\node at (4.2, -0.3) {$1$}; 
\node at (0, -5) {$U_{4,1}$};
\end{tikzpicture}
\qquad 
\begin{tikzpicture}[scale=0.3]
\draw[thick] (2, 3.464)--(4, 0)--(2, -3.464) arc (30:150:2.31) -- (2, 3.464) 
 --cycle;

\draw (0,0) circle(4);
\node at (0, 0) {\Huge $\cdot$}; 
\node at (0, -0.5) {${\rm O}$}; 
\node at (4.2, -0.3) {$1$}; 
\node at (0, -5) {$U_{5,1}$};
\end{tikzpicture}

\end{center}

\caption{$U_{k, 1}$, $1 \le k \le 5$}
\end{figure}
%%%%%%%%%%%%%%%%%%
%%%%%%%%%%%%%%%%%%
%%%%%%%%%%%%%%%%%%%%%%%
%%%%%%%%%%%%%%%%%%%%%%%
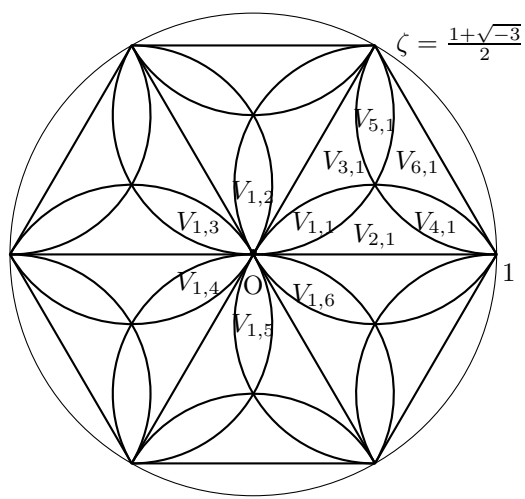
\begin{figure}  \label{fig-3-2}
\begin{center}
\begin{tikzpicture}[scale=0.8]
\draw[thick] (4,0)--(2, 3.464)--(-2, 3.464)--(-4,0)--(-2, -3.464)--(2, -3.464)--(4, 0)--
cycle;
\draw[thick] (2, 3.464) arc (-30:-150:2.31);
\draw[thick] (-2, 3.464) arc (30:-90:2.31);
\draw[thick] (-4, 0) arc (90:-30:2.31);
\draw[thick] (-2, -3.464) arc (150:30:2.31);
\draw[thick] (2, -3.464) arc (210:90:2.31);
\draw[thick] (4, 0) arc (270:150:2.31);

\draw[thick] (4, 0) -- (-4, 0);
\draw[thick] (2, 3.464) -- (-2, -3.464); 
\draw[thick] (-2, 3.464) -- (2, -3.464); 

\draw[thick] (4, 0) arc (-30:-270: 2.31);
\draw[thick](2, -3.464) arc (-30:210:2.31);
\draw[thick] (-4, 0) arc (-150:90:2.31);

\draw[thick] (4, 0) arc (30:270: 2.31);
\draw[thick](2, 3.464) arc (30:-210:2.31);
\draw[thick] (-4, 0) arc (150:-90:2.31);

\draw (0,0) circle(4);
\node at (0, 0) {\Huge $\cdot$}; 
\node at (0, -0.5) {${\rm O}$}; 
\node at (4.2, -0.3) {$1$}; 
\node at (3.4, 3.52) {$\zeta = \frac{1 + \sqrt{-3}}{2}$} ;
\node at (1, 0.5){$V_{1,1}$};
\node at (2, 0.3){$V_{2,1}$};
\node at (1.5, 1.5){$V_{3,1}$};
\node at (3, 0.5){$V_{4,1}$};
\node at (2, 2.2){$V_{5,1}$};
\node at (2.7, 1.5){$V_{6, 1}$};
\node at (0, 1){$V_{1,2}$};
\node at (-0.9, 0.5){$V_{1,3}$};
\node at (-0.9, -0.5){$V_{1,4}$}; 
\node at (0, -1.2){$V_{1, 5}$};
\node at (1, -0.7){$V_{1,6}$};
\end{tikzpicture}
\end{center}
\caption{Finite range structure}
\end{figure}
%%%%%%%%%%%%%%%%%%
%%%%%%%%%%%%%%%%%%
\begin{thm} \label{thm-3-1}
For any cylinder set $\langle a_{1}, a_{2}, \ldots, a_{n} \rangle$ which is not full, 
$T^{n}(\langle a_{1}, a_{2}, \ldots, a_{n} \rangle)$ is one of $U_{k, \ell}$, 
$1 \le k \le 5$, $1 \le \ell \le 6$, a finite union of elements of 
$\mathcal V$ except for a set of Lebesgue measure $0$.  This means 
that $T$ has a finite range structure.  
\end{thm}
{\bf Proof.}    
If $|\alpha|  >  \sqrt{3}$ for $\alpha \in \mathcal J$, then 
$\langle \alpha \rangle$ is full, since $U + \alpha \subset U^{-1}$ in this case, 
see Fig.\ref{fig-1-1}. 
We put $\eta_{k} = \zeta^{(k-1)} \eta$ for $1 \le k \le 6$.  
It is easy to see that $|\alpha|= \sqrt{3}$ means $\alpha = \eta_{k}$ for some 
$1 \le k \le 6$.  
We have 
\begin{equation} \label{eq-3-1}
T(\langle\eta_{k}\rangle) = U_{1, k} \quad \mbox{and} \quad 
T^{2}(\langle \eta_{k}, \eta_{\ell}\rangle) = U_{2, \ell} \,\, \mbox{if }\,\, 
k + \ell = 4 \,\,\mbox{(mod. 6)} .  
\end{equation}
Moreover,  we see, by simple calculation, 
\begin{equation} \label{eq-3-2}
T^{3}(\langle \eta_{2}, \eta_{2}, \eta_{1}+3\rangle) = U_{3, 6}, \,\, 
T^{3}(\langle \eta_{2}, \eta_{2}, \eta_{3}\rangle) = U_{4, 3} \,\,  \mbox{and} \,\, 
T^{3}(\langle \eta_{2}, \eta_{2}, \eta_{1}\rangle) = U_{5, 6} .    
\end{equation}
Also other $U_{k, \ell}$'s for $k = 3, 4, 5$, are  obtained by other cylinder sets of the 
form $\langle \eta_{j_{1}}, \eta_{j_{2}}, \eta_{j_{3}} \rangle$.  Hence we can get all 
$U_{k, \ell}$, $1 \le k \le 5$ and $1 \le \ell \le 6$ as $T^{n}(\langle a_{1}, a_{2}, 
\ldots , a_{n}\rangle)$.  We show that $T^{n}(\langle a_{1}, a_{2}, \ldots , a_{n}
\rangle)$ is one of these $30$ areas or $U$ if it is admissible.    

From the definition of the cylinder set,  we have 
\begin{equation} \label{eq-3-3}
T^{n+1} (\langle a_{1}, a_{2}, \ldots, a_{n}, a_{n+1} \rangle) 
= T\left( T^{n}( \langle a_{1}, a_{2}, \ldots, a_{n} \rangle)\cap \langle a_{n+1} 
\rangle  \right)   
\end{equation} 
for any sequence$(a_{1}, a_{2}, \ldots, a_{n}, a_{n+1})$.  
Hence, $T^{n+1}(\langle a_{1}, a_{2}, \ldots, a_{n}, a_{n+1} \rangle) = 
T(\langle a_{n+1} \rangle)$ whenever $\langle a_{1}, a_{2}, \ldots, a_{n}\rangle$ 
is full.  
Suppose that $(a_{1}, a_{2}, \ldots, a_{n})$ is admissible but not full.  
It is enough to consider the following $6$ cases:  
$T^{n}(\langle a_{1}, a_{2}, \ldots, a_{n} \rangle) = U_{k, 1}$, $1 \le k \le 6$, since 
other cases are obtained by rotations by $\zeta^{\ell}$, $1 \le \ell \le 5$. \\
(i) $k=1$ \\
From \eqref{eq-3-3}, we have 
$T^{n+1}(\langle a_{1}, a_{2}, \ldots, a_{n}, a_{n+1} \rangle) = 
T(U_{1,1} \cap (\langle a_{n+1}\rangle)$.  Thus we have 
\[
T^{n+1}(\langle a_{1}, a_{2}, \ldots, a_{n}, a_{n+1} \rangle) = 
\left\{ \begin{array}{lll} 
U_{2, 3} & \mbox{if} & a_{n+1} = \eta_{3} \\
T(\langle a_{n+1} \rangle) & \mbox{otherwise}
\end{array}\right.  . 
\]
 (ii) $k = 2$ \\
Here our condition means $T^{n+1}(\langle a_{1}, a_{2}, \ldots, a_{n}, a_{n+1} \rangle) = T(U_{2,1} \cap (\langle a_{n+1}\rangle)$. 
Thus we get 
\[
T^{n+1}(\langle a_{1}, a_{2}, \ldots, a_{n}, a_{n+1} \rangle) = 
\left\{ \begin{array}{lll} 
U_{4, 4} & \mbox{if} & a_{n+1} = \eta_{4} \\
U_{5, 1} & \mbox{if} & a_{n+1} = \eta_{2} \\
U_{3,4}  & \mbox{if} & a_{n+1} = \eta_{2}+3j \zeta, \,\,\mbox{or} \,\,
\eta_{4} - 3j \zeta, \,\, j \ge 1 \\
T(\langle a_{n+1} \rangle) & \mbox{otherwise}
\end{array} \right.  .  
\]
\\
(iii) $k=3$ \\
In this case, we have 
\[
T^{n+1}(\langle a_{1}, a_{2}, \ldots, a_{n}, a_{n+1} \rangle) = 
\left\{ \begin{array}{lll} 
U_{3,2} & \mbox{if} & a_{n+1} = \pm 3j \overline{\zeta} \,\,  j \ge 1  \\
T(\langle a_{n+1} \rangle) & \mbox{otherwise}
\end{array} \right. .
\]
\\
(iv) k =4 and 5 \\
It would be sufficient to check the case $k = 4$.  
If $T^{n+1}(\langle a_{1}, a_{2}, \ldots, a_{n}, a_{n+1} \rangle) = T(U_{4,1} \cap (\langle a_{n+1}\rangle)$, then it is easy to check that 
\[
T^{n+1}(\langle a_{1}, a_{2}, \ldots, a_{n}, a_{n+1} \rangle) = 
\left\{ \begin{array}{lll} 
U_{3,2} &   \mbox{if} & a_{n+1} = \pm 3j \overline{\zeta}\,\, j \ge 1 \\
T(\langle a_{n+1} \rangle) & \mbox{otherwise}
\end{array} \right. .
\]
\\
From (i) - (iv) and \eqref{eq-3-1} altogether, we get the conclusion. \\
Note that we only use the following elementary facts: 
\[
\begin{array}{l}
\left\{z : |z- (\pm\tfrac{2}{3} \eta)| = \tfrac{\sqrt{3}}{3}\right\}^{-1}  
=  \left\{z : |z- (\pm \tfrac{2}{3}\overline{\eta})| = \tfrac{\sqrt{3}}{3}\right\} \\
\left\{z : |z - (\tfrac{2}{3}\sqrt{-3})| = \tfrac{\sqrt{3}}{3}\right\}^{-1} = 
\left\{z : |z - (-\tfrac{2}{3}\sqrt{-3})| = \tfrac{\sqrt{3}}{3}\right\}\\
\left\{z : |z- (\pm \tfrac{1}{3} \eta)| = \tfrac{\sqrt{3}}{3}\right\}^{-1} = 
\{z = x + yi : y = \sqrt{3} x \mp \sqrt{3} \} \\
\left\{z : |z- (\pm\tfrac{1}{3} \overline{\eta})| = \tfrac{\sqrt{3}}{3}\right\}^{-1} = 
\{z = x + yi : y = - \sqrt{3} x \pm \sqrt{3} \} \\
\left\{z : |z - (\pm \tfrac{1}{3}\sqrt{-3})| = \tfrac{\sqrt{3}}{3}\right\}^{-1} =
\left\{z = x + yi : y = \mp \tfrac{\sqrt{3}}{2} \right\}\\[8pt]
\{z = x + yi : y =\pm \sqrt{3} x \}^{-1} = \{z = x + yi : y =\mp \sqrt{3} x \}  \\[5pt]
(\mathbb R\cup\{\infty\})^{-1} = \mathbb R \cup \{\infty\},
\end{array}
\]
where $\tfrac{1}{0} = \infty$ and $\tfrac{0}{\infty}$. 
Also note that If $|\alpha| \ge 3 $ for $\alpha \in \mathcal J$, then it is easy to see 
that $\langle \alpha \rangle$ is full, since $U + \alpha \subset U^{-1} = 
\{z : z^{-1} \in U \} $ in this case. It is clear that $U_{k, \ell}$ is a union 
of elements of $\mathcal V$ except for a set of Lebesgue measure $0$. 
\qed \\

In the rest of this section, we discuss the ergodicity of $T$ based on 
M.~Yuri \cite{Y}.  The proof is a sort of routine work and so we only give 
a brief sketch of the proof of the ergodicity.  However, we need the fact 
$|q_{n-1}(z)| <|q_{n}(z)|$ for any admissble sequence $(b_{1}(z), b_{2}(z), \ldots , b_{n}(z))$ which is given in the next section implicitly.  So the following proof will be completed after some discussion in \S4.  We give a comment at the beginning of \S5 concerning this point.   

We rewrite the conditions (c.1), (c.2) and (c.3) in Yuri \cite{Y} for $T$.  We consider 
a cylinder set $\langle a_{1}, \ldots , a_{n} \rangle$ with an admissible sequence 
$(a_{1}, \ldots , a_{n})$.  We put $z_{n} = T^{n}(z)$ and denote the local inverse 
map of $T^{n}$ by $\Psi :z_{n} \mapsto z$.  
A cylinder set $\langle a_{1}, \ldots , a_{n} \rangle$ is said to be a 
$C$-R\'enyi cylinder for a positive constant $C$ if
\[
\sup_{w \in T^{n}\langle a_{1}, \ldots , a_{n} \rangle} |D(\Psi^{n})(w)| 
< C \cdot \inf_{w \in \langle a_{1}, \ldots , a_{n} \rangle} |D(\Psi^{n})(w)| 
\]
holds, where $D\Psi$ denotes the Jacobian of $\Psi$ as a real two dimension 
map.  We denote by $R(C, T)$ the set of $C$-cylinder sets. 
Then the following is our version of conditions:\\
(c.1) (generator condition) The set of cylinder sets separate points, i.e., for any pair $z_{1}$ and $z_{2}\in U$, 
there exists $k\ge 1$ such that $b_{k}(z_{1}) \ne b_{k}(z_{2})$. \\
(c.2) (transitivity condition) For any $1 \le k \le 5$, $1 \le \ell \le 6$, there exists a cylinder set $B =\langle a_{1}, a_{2}, \ldots ,a_{k}\rangle$ of length $k$ such that $B\subset U_{k, \ell}$ and $T^{k}B^{\circ} = U^{\circ}$. \\
(c.3) For any cylinder set $B=\langle a_{1}, a_{2}, \ldots , a_{k}\rangle\in R(C, T)$ 
and $(a_{1}', a_{2}', \ldots, a_{m}')$, $\langle a_{1}', a_{2}', \ldots, a_{m}',    
 a_{1}, a_{2}, \ldots , a_{k}\rangle \in R(C, T)$ whenever it is admissible. 

%%%%%%%%%%%%%%%
%%%%%%%%%%%%%%%
\begin{prop}
The map $T$ is ergodic with respect to the absolutely continuous invariant measure 
$\mu$ where ergodic means that if $T^{-1}A = A$, then $\mu(A) =0$ or 
$\mu(A^{C}) = 0$, where $A^{C}$ denotes the complement of $A$ in $U$.  
%More precisely, $T$ is exact, i. e. If a measurable set $B \subset U$ is in $\bigcap_{n=1}^{\infty} T^{-n}{\mathbb B}$, then $\mu(B) =0$ or $\mu(B^{C})=0$. 
\end{prop}
%%%%%%%%%%%%%%%
%%%%%%%%%%%%%%%
{\bf Proof.}  
First we show that the conditions (c.1), (c.2) and (c.3) hold for $T$ after a simple 
modification.    
We can choose 
\[
\tilde{R}(T:3) = \{\langle a_{1}, \ldots , a_{n} \rangle : n \ge 1, \,\, 
(a_{1}, \ldots , a_{n} )\,\,{\rm is \,\,admissible},\,\, |a_{n}| \ge 3 \}
\]
with $C=9$.  This is possible since 
\[
\Psi (z_{n}) = \frac{p_{n-1} z_{n} + p_{n}}{q_{n-1} z_{n} + q_{n}}
\]
for $z \in \langle a_{1}, \ldots , a_{n} \rangle$ with $z_{n} = T^{n}(z)$ and then 
\[
(\Psi)'(z_{n}) = \frac{1}{(q_{n} + q_{n-1} z_{n})^{2}} .
\]
As we tentatively assume that $|q_{n-1}(z)| > |q_{n-2}(z)|$ holds for any admissible 
sequence $(a_{1}, \ldots , a_{n-1} )$,  we see 
\[
|q_{n}(z)|> 2|q_{n-1}(z)| \qquad \mbox{if} \,\, |a_{n}|\ge 3 . 
\] 
Hence we see 
\[
\frac{1}{3 |q_{n}|^{2}} < \left|(\Psi)'(z_{n})\right| < \frac{1}{|q_{n}|^{2}} .  
\]
Since $|D\Psi| = |(\Psi)'|^{2}$, we get 
\[
\tilde{R}(T:3) \subset R(9, T) .  
\]
In the proof in \cite{Y}, it is easy to see that one can replace $R(9, T)$ by $\tilde{R}(T:3)$.   
Then the conditions (c.1), (c.2) and (c.3) are easy to follow.  Indeed, 
(c.1) is equivalent to the convergence of the continued fraction expansions, 
(c.2) is also easily checked, and (c.3) is trivial.  

Once we have these three properties, we can define a jump map as 
\[
T_{J}(z) = \left\{ \begin{array}{lll}
T^{N_{J}+1}(z) & \mbox{if}& N_{J} <\infty \\
z & \mbox{if}& N_{J} =\infty \end{array} \right.  
\]
with $N_{J}=N_{J}(z) = \min\{ n\ge 1 : |b_{n}|\ge 3\}$ if minimum $n$ exists, otherwise 
$N_{J}(z) = \infty$.   Then, from \cite{I-Y},  we have the existence of the 
absolutely continuous invariant measure $\mu_{J}$ for $T_{J}$.  Moreover, we 
can induce the absolutely continuous invariant measure for $T$, which may not be a finite measure.  The condition (c.2) implies that the support of the invariant 
measure is $U$ and that the ergodicity of $T$.   
\qed\\
  
\begin{rem} 
(i) The condition (c.4) in \cite{Y} is the necessary and sufficient condition 
for the absolutely continuous invaroant measure to be finite.  In the next section, 
we will find that $\mu$ is a finite measure.  Thus, once we get the finiteness, the condition (c.4) automatically follows.  Then because of \cite[Theorem 1]{Y}, 
the map is exact.  \\
(ii) If we follow \cite{Y} to get the Rokholin's formula for the entropy $h(T)$: 
\[
\int_{U} \log |DT| d\mu, 
\]
then we need to check some more conditions as well as the weak Bernoulli property. 
However, it is a sort of routine work as we mentioned before and so 
we show the existence of the L\'evy's constant without Rokholin's formula in the 
last part of the next section.
\end{rem}

%%%%%%%%%%%%%%%%%%%%%%%%
%%%%%%%%%%%%%%%%%%%%%%%%
%%%%%%%%%%%%%%%%%%%%%%%%
%%%%%%%%%%%%%%%%%%%%%%%%

\section{Dual area}
\begin{figure} \label{fig-4-1}
\begin{center}
\begin{tikzpicture}[scale=0.3]
\fill[lightgray] (-2, 3.464) arc (120:300:4) arc (-120:60:2) arc (-60:120:2) arc 
(0:180:2);

\draw (4,0)--(2, 3.464)--(-2, 3.464)--(-4,0)--(-2, -3.464)--(2, -3.464)--(4, 0)--
cycle;
\draw (-2, 3.464) arc (120:300:4);
\draw (4, 0) arc (-60:120:2);
\draw (4, 0) arc (60:-120:2);
\draw (2, 3.464) arc (0:180:2);
\node at (0, 0) {\Large $\cdot$}; 
\node at (0, -0.6) {\bf 0};

\node at (0, -5.5) {$V_{0, 1,1}^{*}$}; 
\end{tikzpicture}
%%%%%%%%%%%%%%%%%%%%
%%%%%%%%%%%%%%%%%%%%
\qquad
\begin{tikzpicture}[scale=0.3]
\fill[lightgray] (4, 0) arc (-60:120:2) arc (60:300:4) arc (-120:60:2);

\draw (4,0)--(2, 3.464)--(-2, 3.464)--(-4,0)--(-2, -3.464)--(2, -3.464)--(4, 0)--
cycle;
\draw (2, 3.464) arc (60:300:4);
\draw (4, 0) arc (-60:120:2);
\draw (4, 0) arc (60:-120:2);
\node at (0, 0) {\Large $\cdot$}; 
\node at (0, -0.6) {\bf 0};

\node at (0, -5.5) {$V_{0, 2,1}^{*}$}; 
\end{tikzpicture}
%%%%%%%%%%%%%%%%%%%%
%%%%%%%%%%%%%%%%%%%%
\qquad
\begin{tikzpicture}[scale=0.3]
\fill[lightgray] (4, 0) arc (-60:120:2) arc (0:180:2) arc (120:360:4);

\draw (4,0)--(2, 3.464)--(-2, 3.464)--(-4,0)--(-2, -3.464)--(2, -3.464)--(4, 0)--
cycle;
\draw (-2, 3.464) arc (120:360:4);
\draw (4, 0) arc (-60:120:2);
\draw (2, 3.464) arc (0:180:2);
\node at (0, 0) {\Large $\cdot$}; 
\node at (0, -0.6) {\bf 0};

\node at (0, -5.5) {$V_{0, 3,1}^{*}$}; 
\end{tikzpicture} 
\\
%%%%%%%%%%%%%%%%%%%%
%%%%%%%%%%%%%%%%%%%%
\begin{tikzpicture}[scale=0.3]
\fill[lightgray] (4, 0) arc (-120: 180:4) arc (60:300:4) arc (-120:60:2);

\draw (4,0)--(2, 3.464)--(-2, 3.464)--(-4,0)--(-2, -3.464)--(2, -3.464)--(4, 0)--
cycle;
\draw (2, 3.464) arc (60:300:4);
\draw (4, 0) arc (-120:180:4);
\draw (4, 0) arc (60:-120:2);
\node at (0, 0) {\Large $\cdot$}; 
\node at (0, -0.6) {\bf 0};

\node at (0, -5.5) {$V_{0, 4,1}^{*}$}; 
\end{tikzpicture}
%%%%%%%%%%%%%%%%%%%%
%%%%%%%%%%%%%%%%%%%%
%%%%%%%%%%%%%%%%%%%%
\qquad
\begin{tikzpicture}[scale=0.3]
\fill[lightgray] (-2, 3.464) arc (120:360:4) arc (-120:180:4) arc (0:180:2);

\draw (4,0)--(2, 3.464)--(-2, 3.464)--(-4,0)--(-2, -3.464)--(2, -3.464)--(4, 0)--
cycle;
\draw (-2, 3.464) arc (120:360:4);
\draw (4, 0) arc (-120:180:4);
\draw (2, 3.464) arc (0:180:2);
\node at (0, 0) {\Large $\cdot$}; 
\node at (0, -0.6) {\bf 0};

\node at (0, -5.5) {$V_{0, 5,1}^{*}$}; 
\end{tikzpicture}
%%%%%%%%%%%%%%%%%%%%
%%%%%%%%%%%%%%%%%%%%
\quad
\begin{tikzpicture}[scale=0.3]
\fill[lightgray] (2, 3.464) arc (60:360:4) arc (-120:180:4);
\draw (4,0)--(2, 3.464)--(-2, 3.464)--(-4,0)--(-2, -3.464)--(2, -3.464)--(4, 0)--
cycle;
\draw (2, 3.464) arc (60:360:4);
\draw (4, 0) arc (-120:180:4);
\node at (0, 0) {\Large $\cdot$}; 
\node at (0, -0.6) {\bf 0};

\node at (0, -5.5) {$V_{0, 6,1}^{*}$}; 
\end{tikzpicture}
\end{center}
\caption{$V_{0, k, 1}^{*}$, $1 \le k \le 6$ : outside of the gray zone}

\end{figure}
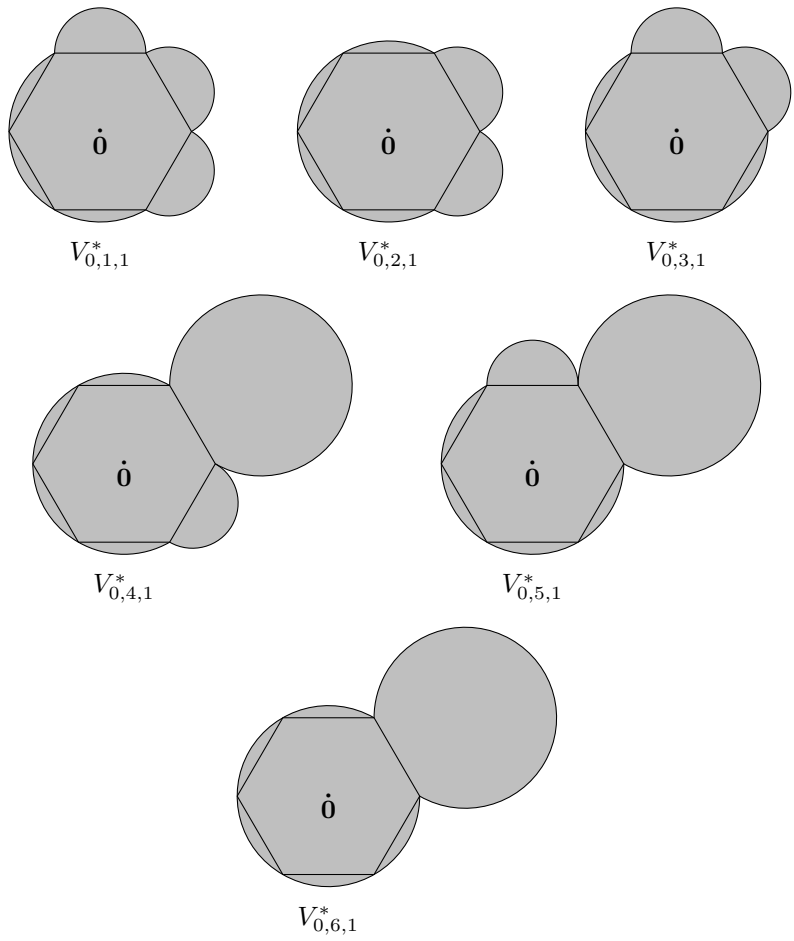

Now we consider the ``dual sequences" $\left(-\tfrac{q_{n-1}}{q_{n}}: n \ge 1\right)$ 
of the convergents $\tfrac{p_{n}}{q_{n}}$.  The same as in 
\cite[page 50--51]{E-I-N-N}, we have for any $(k, \ell)$, $1 \le k, \ell \le 6$, 
\[
V_{k, \ell}^{*} := \{-\tfrac{q_{n}(z)}{q_{n-1}(z)} : n \ge 1, \,\, T^{n}(z) = z_{0} \in 
V_{k, \ell} \}^{cl}
\]
is independent of the choice of $z_{0} \in V_{k, \ell}$, $1 \le k, \ell \le 6$, where 
$A^{cl}$ denotes the closure of $A \subset \mathbb C$. 
In the following, we will give $V_{k, \ell}^{*}$, $1 \le k, \ell \le 6$ explicitly, which allows us to see the monotonicity of $(|q_{n}(z)| : n \ge 1)$ for any $z \in U$, see \S5, and also to show some ergodic properties of $T$.    
%%%%%%%%%%%%%%%%%%%%%%%%%%
%%%%%%%%%%%%%%%%%%%%%%%%%%
\begin{defn} We define 
\[
\begin{array}{lll}
V_{0,1, 1} ^{*}&=& \{z \in \mathbb C: |z|>1, 
                         \left|z - \tfrac{\sqrt{-3}}{2}\right| > \tfrac{1}{2},  
                         \left|z - \tfrac{\eta}{2}\right|>\tfrac{1}{2}, 
                         \left|z - \tfrac{\overline{\eta}}{2}\right|>\tfrac{1}{2} \} \\
V_{0,2, 1} ^{*}&=& \{z \in \mathbb C: |z|>1, 
                          \left|z - \tfrac{\eta}{2}\right|>\tfrac{1}{2}, 
                         \left|z - \tfrac{\overline{\eta}}{2}\right|>\tfrac{1}{2} \} \\
V_{0, 3, 1} ^{*}&=& \{z \in \mathbb C: |z|>1, 
                         \left|z - \tfrac{\sqrt{-3}}{2}\right| > \tfrac{1}{2},  
                         \left|z - \tfrac{\eta}{2}\right|>\tfrac{1}{2} \} \\ 
 V_{0,4, 1} ^{*}&=& \{z \in \mathbb C: |z|>1, 
                          \left|z - \eta\right|> 1, 
                         \left|z - \tfrac{\overline{\eta}}{2}\right|>\tfrac{1}{2} \} \\
 V_{0, 5, 1} ^{*}&=& \{z \in \mathbb C: |z|>1, 
                          \left|z - \eta\right|> 1, 
                         \left|z - \tfrac{\sqrt{-3}}{2}\right|>\tfrac{1}{2} \} \\
 V_{0,6, 1} ^{*}&=& \{z \in \mathbb C: |z|>1, 
                          \left|z - \eta\right|> 1\}
\end{array}.
\]
and 
\[
V_{0, k, \ell}^{*} = \zeta^{(\ell -1)}V_{0, k, 1}^{*}, \,\,1 \le k \le 6, \,\, 
1\le \ell \le 6 ,  
\]
see {\rm Fig. \ref{fig-4-1}}.
\end{defn}

%%%%%%%%%%%%%%%%%%%%
%%%%%%%%%%%%%%%%%%%%
%%%%%%%%%%%%%%%%%%%%
%%%%%%%%%%%%%%%%%%%%
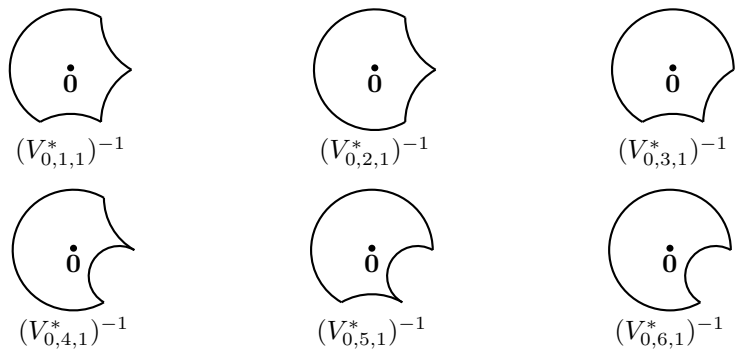
\begin{figure} \label{fig-4-2}
\begin{center}
\begin{tikzpicture}[scale=0.2]
\draw[thick, black] (2, 3.464) arc (60:240:4);
\draw[thick, black] (4, 0) arc (240:180:4);
\draw[thick, black] (4, 0) arc (120:180:4);
\draw[thick, black] (2, -3.464) arc (60:120:4);

\node at (0, 0) {\Huge $\cdot$};
\node at (0, -0.9) {\bf 0};
\node at (0, -5.5) {$(V_{0,1,1}^{*})^{-1}$};

\end{tikzpicture}
\qquad \qquad \qquad
\begin{tikzpicture}[scale=0.2]

\draw[thick, black] (2, 3.464) arc (60:300:4);
\draw[thick, black] (4, 0) arc (240:180:4);
\draw[thick, black] (4, 0) arc (120:180:4);
\node at (0, 0) {\Huge $\cdot$};
\node at (0, -0.9) {\bf 0};
\node at (0, -5.5) {$(V_{0,2,1}^{*})^{-1}$};
\end{tikzpicture}
%%%%%%%%%%%%%%%%%%%%
\qquad\qquad \qquad
%%%%%
%%%%%%%%%%%%%%%%%%%%
\begin{tikzpicture} [scale=0.2]
\draw[thick, black] (4, 0) arc (0:240:4);

\draw[thick, black] (4, 0) arc (120:180:4);
\draw[thick, black] (2, -3.464) arc (60:120:4);

\node at (0, 0) {\Huge $\cdot$};
\node at (0, -0.9) {\bf 0};
\node at (0, -5.5) {$(V_{0,3,1}^{*})^{-1}$};

\end{tikzpicture}
\qquad \qquad \qquad
\begin{tikzpicture} [scale=0.2]

\draw[thick, black] (2, 3.464) arc (60:300:4);
\draw[thick, black] (4, 0) arc (240:180:4);

\draw[thick, black] (4, 0) arc (60:240:2);
\node at (0, 0) {\Huge $\cdot$};
\node at (0, -0.9) {\bf 0};
\node at (0, -5.5) {$(V_{0,4,1}^{*})^{-1}$};
\end{tikzpicture}
\qquad\qquad \qquad
\begin{tikzpicture} [scale=0.2]
\draw[thick, black] (4, 0) arc (0:240:4);

\draw[thick, black] (2, -3.464) arc (60:120:4);

\draw[thick, black] (4, 0) arc (60:240:2);

\node at (0, 0) {\Huge $\cdot$};
\node at (0, -0.9) {\bf 0};
\node at (0, -5.5) {$(V_{0,5,1}^{*})^{-1}$};

\end{tikzpicture}
\qquad \qquad \qquad
\begin{tikzpicture} [scale=0.2]

\draw[thick, black] (4, 0) arc (0:300:4);

\draw[thick, black] (4, 0) arc (60:240:2);

\node at (0, 0) {\Huge $\cdot$};
\node at (0, -0.9) {\bf 0};
\node at (0, -5.5) {$(V_{0,6,1}^{*})^{-1}$};

\end{tikzpicture}
\end{center}
\caption{$(V_{0, k, 1}^{*})^{-1}$, $1 \le k \le 6$}

\end{figure}

%%%%%%%%%%%%%%%%%%%%
%%%%%%%%%%%%%%%%%%%%
We define a map $\hat{T}$ on $U \times (\mathbb C \cup \{\infty \})$ by 
\begin{equation} \label{eq-4-1}
\hat{T}(z, w) = \left(\tfrac{1}{z} - b_{1}(z), \tfrac{1}{w} - b_{1}(z) \right) \,\, 
\mbox{for} \,\, (z, w) \in U \times (\mathbb C \cup \{\infty\}).  
\end{equation}
Then it is easy to see the following. 
\begin{lem} \label{lem-4-1}
For any $z \in U$, we have 
\[
\hat{T}^{n}(z, -\infty) = (T^{n}(z), -\tfrac{q_{n}(z)}{q_{n-1}(z)}) 
\]
for $n \ge 0$, where $\tfrac{q_{0}(z)}{q_{-1}(z)} = \tfrac{1}{0} = \infty$.

\end{lem}
{\bf Proof.}  It follows from the definition of $\hat{T}$, 
\[
\hat{T}(z, \infty) = (T(z), -b_{1}(z)) = \left(T(z), -\tfrac{q_{1}(z)}{q_{0}(z)}\right). 
\]
Then we get 
\[
\hat{T}\left(T^{n}(z), -\tfrac{q_{n}(z)}{q_{n-1}(z)}\right) = 
\left(T^{n+1}(z), -\tfrac{q_{n-1}(z)}{q_{n}(z)} - b_{n+1}(z)\right) =  
\left(T^{n+1}(z), 
-\tfrac{q_{n+1}(z)}{q_{n}(z)}\right)
\]
for $n \ge 1$ by induction.  
\qed
%%%%%%%%%%%%%%%%%%%%%%%
%%%%%%%%%%%%%%%%%%%%%%%
%%%%%%%%%%%%%%%%%%%%%%%
%%%%%%%%%%%%%%%%%%%%%%%

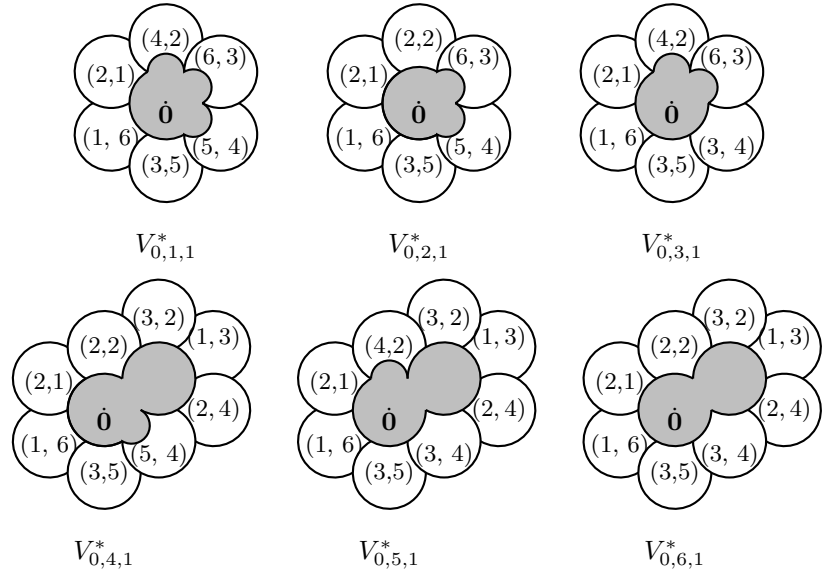
\begin{figure} \label{fig-4-3}
\begin{center}
\begin{tikzpicture}[scale=0.12]
\fill[lightgray] (-2, 3.464) arc (120:300:4) arc (-120:60:2) arc (-60:120:2)  
arc (0:180:2);

\draw[thick] (4, 0) arc (-60:120:2) arc (180:-120: 4);
\draw[thick] (4, 0) arc (60:-120:2) arc (-180: 60:4); 
\draw[thick] (-2, 3.464) arc (120:300:4);
\draw[thick] (-2, -3.464) arc (-240:0:4);
\draw[thick] (-4, 0) arc (-60:-120:4) arc (120:300:4)   ;
\draw[thick] (-8, 0) arc (240:60:4);
\draw[thick] (-2, 3.464) arc (240:0:4);
\draw[thick] (2, 3.464) arc (0:180:2);

\node at (0, 7) {\small (4,2)};
\node at (-6, 3) {\small (2,1)};
\node at (-6, -4) {\small (1, 6)};
\node at (0, -7) {\small (3,5)};
\node at (6, -5) {\small (5, 4)};
\node at (6,5) {\small $(6, 3)$}; 

\node at (0, 0) {$\cdot$};
\node at (0,-1.3) {\bf 0};
\node at (0, -16) {$V _{0, 1, 1}^{*}$};
\end{tikzpicture}
\qquad
%%%%%%%%%%%%%%%%%%%
%%%%%%%%%%%%%%%%%%%%
\begin{tikzpicture}[scale=0.12]
\fill[lightgray] (2, 3.464) arc (60:300:4) arc (-120: 60: 2) arc (-60: 120:2);

\draw[thick] (4, 0) arc (-60:120:2) arc (180:-120: 4);
\draw[thick] (4, 0) arc (60:-120:2) arc (-180: 60:4); 
\draw[thick] (-2, 3.464) arc (120:300:4);
\draw[thick] (-2, -3.464) arc (-240:0:4);
\draw[thick] (-4, 0) arc (-60:-120:4) arc (120:300:4)   ;
\draw[thick] (-8, 0) arc (240:60:4);
\draw[thick] (-2, 3.464) arc (240:0:4);

\draw[thick] (2, 3.464) arc (60:300:4);

\node at (0, 7) {\small (2,2)};
\node at (-6, 3) {\small (2,1)};
\node at (-6, -4) {\small (1, 6)};
\node at (0, -7) {\small (3,5)};
\node at (6, -5) {\small (5, 4)};
\node at (6,5) {\small $(6, 3)$}; 

\node at (0, 0) {$\cdot$};
\node at (0,-1.3) {\bf 0};
\node at (0, -16) {$V_{0, 2, 1}^{*}$};
\end{tikzpicture}
%%%%%%%%%%%%%%%%%%%
%%%%%%%%%%%%%%%%%%%%
\qquad
\begin{tikzpicture}[scale=0.12]

\fill[lightgray] (-2, 3.464) arc (120:360:4) arc (-60:120:2) 
arc (0:180:2);

\draw[thick] (4, 0) arc (-60:120:2) arc (180:-120: 4);

\draw[thick] (2, -3.464) arc (-180: 60:4); 
\draw[thick] (-2, 3.464) arc (120:300:4);
\draw[thick] (-2, -3.464) arc (-240:0:4);
\draw[thick] (-4, 0) arc (-60:-120:4) arc (120:300:4)   ;
\draw[thick] (-8, 0) arc (240:60:4);
\draw[thick] (-2, 3.464) arc (240:0:4);

\draw[thick] (4, 0) arc (0:-60:4);
\draw[thick] (2, 3.464) arc (0:180:2);

\node at (0, 7) {\small (4,2)};
\node at (-6, 3) {\small (2,1)};
\node at (-6, -4) {\small (1, 6)};
\node at (0, -7) {\small (3,5)};
\node at (6, -5) {\small (3, 4)};
\node at (6,5) {\small $(6, 3)$};

\node at (0, 0) {$\cdot$};
\node at (0,-1.3) {\bf 0};
\node at (0, -16) {$V_{0, 3, 1}^{*}$};

\end{tikzpicture}
\begin{tikzpicture}[scale=0.12]
\fill[lightgray] (2, 3.464) arc (60:300:4) arc (-120:60:2) arc (-120:180:4);

\draw[thick] (4, 0) arc (-120:180: 4);
\draw[thick] (4, 0) arc (60:-120:2) arc (-180: 60:4); 
\draw[thick] (-2, 3.464) arc (120:300:4);
\draw[thick] (-2, -3.464) arc (-240:0:4);
\draw[thick] (-4, 0) arc (-60:-120:4) arc (120:300:4)   ;
\draw[thick] (-8, 0) arc (240:60:4);
\draw[thick] (-2, 3.464) arc (240:0:4);

\draw[thick] (2, 3.464) arc (60:300:4);
\draw[thick] (10, -3.464) arc (-120:120:4);
\draw[thick] (8, 6.93) arc (-60:180:4);
\draw[thick] (14, 3.464) arc (-60:120:4);

\node at (0, 7) {\small (2,2)};
\node at (-6, 3) {\small (2,1)};
\node at (-6, -4) {\small (1, 6)};
\node at (0, -7) {\small (3,5)};
\node at (6, -5) {\small (5, 4)};
\node at (12, 0) {\small $(2, 4)$}; 
\node at (12, 8) {\small $(1, 3)$}; 
\node at (6, 10) {\small $(3, 2)$}; 

\node at (0, 0) {\small $\cdot$};
\node at (0,-1.3) { \bf 0};
\node at (0, -16) {$V_{0, 4, 1}^{*}$};

\end{tikzpicture}
%%%%%%%%%%%%%%%%%%%
\quad
\begin{tikzpicture}[scale=0.12]
\fill[lightgray] (-2, 3.464) arc (120:360:4) arc (-120:180:4) arc (0:180:2);

\draw[thick] (4, 0) arc (0:-60:4);
\draw[thick] (2, -3.464) arc (-180:60:4);
\draw[thick] (-2, 3.464) arc (120:300:4);
\draw[thick] (-2, -3.464) arc (-240:0:4);
\draw[thick] (-4, 0) arc (-60:-120:4) arc (120:300:4)   ;
\draw[thick] (-8, 0) arc (240:60:4);
\draw[thick] (-2, 3.464) arc (240:0:4);
\draw[thick] (4, 0) arc (-120:180: 4);

\draw[thick] (10, -3.464) arc (-120:120:4);
\draw[thick] (8, 6.93) arc (-60:180:4);
\draw[thick] (14, 3.464) arc (-60:120:4);

\draw[thick] (2, 3.464) arc (0:180:2);

\node at (0, 7) {\small (4,2)};
\node at (-6, 3) {\small (2,1)};
\node at (-6, -4) {\small (1, 6)};
\node at (0, -7) {\small (3,5)};
\node at (6, -5) {\small (3, 4)};
\node at (12, 0) {\small $(2, 4)$}; 
\node at (12, 8) {\small $(1, 3)$}; 
\node at (6, 10) {\small $(3, 2)$}; 

\node at (0, 0) {\small$\cdot$};
\node at (0,-1.3) {\bf 0};
\node at (0, -16) {$V_{0, 5, 1}^{*}$};
\end{tikzpicture}
\quad
%%%%%%%%%%%%%%%%%%%
%%%%%%%%%%%%%%%%%%%%
\begin{tikzpicture}[scale=0.12]
\fill[lightgray] (2, 3.464) arc (60:360:4) arc (-120:180:4);

\draw[thick] (4, 0) arc (0:-60:4);
\draw[thick] (2, -3.464) arc (-180:60:4);   
\draw[thick] (-2, 3.464) arc (120:300:4);
\draw[thick] (-2, -3.464) arc (-240:0:4);
\draw[thick] (-4, 0) arc (-60:-120:4) arc (120:300:4)   ;
\draw[thick] (-8, 0) arc (240:60:4);
\draw[thick] (-2, 3.464) arc (240:0:4);

\draw[thick] (10, -3.464) arc (-120:120:4);
\draw[thick] (8, 6.93) arc (-60:180:4);
\draw[thick] (14, 3.464) arc (-60:120:4);

\draw[thick] (4, 0) arc (-120:180: 4);

\draw[thick] (2, 3.464) arc (60:120:4);

\node at (0, 7) {\small (2,2)};
\node at (-6, 3) {\small (2,1)};
\node at (-6, -4) {\small (1, 6)};
\node at (0, -7) {\small (3,5)};
\node at (6, -5) {\small (3, 4)};
\node at (12, 0) {\small $(2, 4)$}; 
\node at (12, 8) {\small $(1, 3)$}; 
\node at (6, 10) {\small $(3, 2)$}; 

\node at (0, 0) {\small$\cdot$};
\node at (0,-1.3) {\bf 0};
\node at (0, -16) {$V_{0, 6, 1}^{*}$};

\end{tikzpicture}
%%%%%%%%%%%%%%%%%%
%%%%%%%%%%%%%%%%%%

\end{center}
\caption{Configurations around the unit circle: 
$(k, \ell)$ means $(V_{0,k,\ell}^{*})^{-1} - \alpha$ for a suitable $\alpha$}
\end{figure}

We want to show that $V_{k, \ell}^{*} = V_{0, k, \ell}^{*}$ for any $1\le k \le6$ and 
$1\le \ell \le 6$.  First, we show that   $V_{k, \ell}^{*} \subset V_{0, k, \ell}^{*}$.  
For this purpose, we consider their inverse sets $(V_{0,k,\ell}^{*})^{-1}$.  Since 
$\inf_{z \in V_{0, k, \ell}^{*} } |z| = 1$ for any $(k, \ell)$, $1 \le k, \ell \le 6$, 
$\sup_{z \in (V_{0, k, \ell}^{*})^{-1} } |z| = 1$.  Indeed, we have the following. 
%%%%%%%%%%%%%%%%%%
%%%%%%%%%%%%%%%%%%
\begin{lem} \label{lem-4-2} 
The restriction of the map $\hat{T}$ to $\hat{U}= \bigcup_{1 \le k, \ell \le 6} 
V_{k, \ell} \times V_{0, k, \ell}^{*}$ is bijective except for a set of 
(complex-)Lebesgue measure $0$.      
\end{lem}
%%%%%%%%%%%%%
%%%%%%%%%%%%%
{\bf Proof.}  Suppose that $z \in V_{k, \ell}$ and $b_{1}(z)= \alpha$.  Then we have 
$\hat{T}(z, V_{0, k, \ell}^{*})  = (T(z), (V_{0, k, \ell}^{*})^{-1} - \alpha)$.    From the definition of $V_{0, k, \ell}$, it follows that $(V_{0, k, \ell}^{*})^{-1} \subset 
\{z : |z|<1\}$. 
%We consider $V_{0, k, 1}$, $1 \le k\le 6$.  
For $z \in U$, if $b_1(z) = \zeta^m \cdot \eta$, $0 \le m \le 5$, then 
\[
T(z) \notin \bigcup_{k=4}^{6}V_{k, j}, 
\]
where
\[
j = 
\left\{ \begin{array}{lll} 
m + 4 & \mbox{if} & 0 \le m \le 2 \\
m - 2 & \mbox{if} & 3 \le m \le 5
\end{array}\right.  . 
\]
Thus, for each $4 \le k \le 6$, $1 \le \ell \le 6$, there is no 
$k'$ and $\ell'$ such that $(V_{0, k', \ell'}^{*})^{-1} - \zeta^{-(\ell -1)} \eta \subset 
V_{0, k, \ell}^{*}$, see Fig.$6$.  \\
One can check that 
\begin{equation} \label{eq-4-2}
\begin{array}{l}
\left((V_{0, 6, 3}^{*})^{-1}  - \eta_{4}\right)  \cup \left((V_{0, 4, 2}^{*})^{-1}  - 
\eta_{5}\right) \cup 
\left((V_{0, 2, 1}^{*})^{-1}  - \eta_{6} \right)\cup  \\
\left((V_{0, 1, 6}^{*})^{-1}  - \eta_{1} \right)  \cup
\left((V_{0, 3, 5}^{*})^{-1}  - \eta_{2}\right) \cup 
\left((V_{0, 5, 4}^{*})^{-1}  - \eta_{3}  \right)\\ 
\subset  V_{0, 1, 1}^{*} \\ 
{}\\
\left((V_{0, 6, 3}^{*})^{-1}  - \eta_{4} \right) \cup  
\left((V_{0, 2, 2}^{*})^{-1}  - \eta_{5} \right) \cup  
\left((V_{0, 2, 1}^{*})^{-1}  - \eta_{6} \right) \cup \\ 
\left((V_{0, 1, 6}^{*})^{-1}  - \eta_{1}  \right) \cup  
\left((V_{0, 3, 5}^{*})^{-1}  - \eta_{2} \right) \cup 
\left((V_{0, 5, 4}^{*})^{-1}  - \eta_{3}  \right)\\
\subset  V_{0, 2, 1}^{*} \\ 
{}\\
\left((V_{0, 6, 3}^{*})^{-1}  - \eta_{4} \right) \cup 
\left((V_{0, 4, 2}^{*})^{-1}  - \eta_{5} \right)\cup 
\left((V_{0, 2, 1}^{*})^{-1}  - \eta_{6}\right) \cup \\
\left((V_{0, 1, 6}^{*})^{-1}  - \eta_{1}\right) \cup 
\left((V_{0, 3, 5}^{*})^{-1}  - \eta_{2}\right) \cup 
\left((V_{0, 3, 4}^{*})^{-1}  - \eta_{3}\right) \\ 
\subset  V_{0, 3, 1}^{*} \\
{} \\
\left((V_{0,2, 2}^{*})^{-1}  - \eta_{5}\right)\cup 
\left((V_{0, 2, 1}^{*})^{-1}  - \eta_{6} \right)\cup 
\left((V_{0, 1, 6}^{*})^{-1}  - \eta_{1} \right)\cup  \\ 
\left((V_{0, 3, 5}^{*})^{-1}  - \eta_{2} \right)\cup 
\left((V_{0, 5, 4}^{*})^{-1}  - \eta_{3}\right) \cup \\  
\left((V_{0, 2, 4}^{*})^{-1}  - (-3)  \right)\cup 
\left((V_{0, 1, 3}^{*})^{-1}  - (-2 \eta)  \right)\cup
\left((V_{0, 3, 2}^{*})^{-1}  - (-3\zeta)\right)\\
\subset  V_{0, 4, 1}^{*}  
{} \\
\left((V_{0,4, 2}^{*})^{-1}  - \eta_{5}  \right)\cup 
\left((V_{0, 2, 1}^{*})^{-1}  - \eta_{6} \right)\cup 
\left((V_{0, 1, 6}^{*})^{-1}  - \eta_{1} \right)\cup \\
\left((V_{0, 3, 5}^{*})^{-1}  - \eta_{2}  \right)\cup 
\left((V_{0, 3, 4}^{*})^{-1}  - \eta_{3} \right)\cup  \\
\left((V_{0, 2, 4}^{*})^{-1}  - (-3)  \right)\cup
\left((V_{0, 1, 3}^{*})^{-1}  - (-2 \eta) \right) \cup
\left((V_{0, 3, 2}^{*})^{-1}  - (-3\zeta)\right)\\
\subset  V_{0, 5, 1}^{*} \\ 
{}\\
\left((V_{0,2, 2}^{*})^{-1}  - \eta_{5} \right)\cup 
\left((V_{0, 2, 1}^{*})^{-1}  - \eta_{6} \right)\cup 
\left((V_{0, 1, 6}^{*})^{-1}  - \eta_{1} \right)\cup \\
\left((V_{0, 3, 5}^{*})^{-1}  - \eta_{2} \right)\cup 
\left((V_{0, 3, 4}^{*})^{-1}  - \eta_{3} \right)\cup\\
\left((V_{0, 2, 4}^{*})^{-1}  - (-3)\right)  \cup
\left((V_{0, 1, 3}^{*})^{-1}  - (-2 \eta)\right)  \cup
\left((V_{0, 3, 2}^{*})^{-1}  - (-3\zeta)\right)\\
\subset  V_{0, 6, 1}^{*} 
\end{array} , 
\end{equation}
see Fig.\ref{fig-4-3}.  
For other $\alpha \in \eta \cdot \mathfrak o (\sqrt{-3})\setminus{\{0\}}$, it is clear that 
$(V_{0, k, \ell}^{*})^{-1} + \alpha \subset V_{0, k', \ell'}^{*}$, since $|z|\le 1$ for 
$z \in V_{0, k, \ell}$.  By rotations $\zeta^{\ell}$, $0 \le \ell \le 5$.  The same holds for every $(k, \ell)$ thus we see that $\hat{T}$ on $\hat{U}$ is into map.  In 
\eqref{eq-4-2}, it is easy to check that all union are disjoint except for a set of 
boundary points of $(V_{0, k, \ell}^{*})^{-1} - \alpha$, where $\alpha$ is suitably chosen.  It is also not hard to check that the same disjointness hold for all pair of 
$(V_{0, k', \ell'}^{*})^{-1} - \alpha$ and $(V_{0, k'', \ell''}^{*})^{-1} - \alpha'$ in the 
right side of the same $V_{k, \ell}^{*}$ unless 
they are fit to each other, see Fig.$7$.  Here we note that $\zeta (V_{0,k,\ell}^{*})^{-1} 
= (V_{0, k-1, \ell}^{*})^{-1}$, see Fig.$5$. In this way,  we get the assertion of this lemma.
%%%%%%%%%%%%%%%
%%%%%%%%%%%%%%%
%%%%%%%%%%%%%%%
%%%%%%%%%%%%%%%
\begin{figure} \label{fig-4-4}
\begin{center}
\begin{tikzpicture}[scale=0.12]
\fill[lightgray] (4, 6.928) arc (0: 180:4) arc (60: 240:4) arc (120:300:4) arc (180:360:4) 
arc(240:420:4) arc (-60:120:4); 
\fill[gray] (-2, 3.464) arc (120:300:4) arc (-120:60:2) arc (-60:120:2)  
arc (0:180:2);

\draw[thick] (2, 10.392) arc (180:120:4) arc (0:180:4) arc (60:0:4); 
\node at (0, 14) {\small $(1,2)$};

\draw[thick] (-2, 10.392) arc (0: 240:4) ; 
\node at (-6, 11) {\small $(3,1)$}; 

\draw[thick] (-8, 6.928) arc (240:180:4) arc (60:240:4) arc (120:60:4); 
\node at (-12, 7) {\small $(1,1)$};

\draw[thick] (-14, 3.464) arc (120:300:4) ; 
\node at (-12, 0) {\small $(3,6)$};

\draw[thick] (-14, -3.464) arc (120:300:4) arc (180:120:4)   ;
\node at (-12, -7) {\small $(1,6)$};

\draw[thick] (-8, -6.928) arc (120:360:4);  
\node at (-6, -11) {\small $(2, 6)$};

\draw[thick] (-4, -13.856) arc (180:360:4) arc (240:180:4);  
\node at (0, -15) {\small $(1,5)$};

\draw[thick] (2, -10.392) arc (-180: 60:4) ;
\node at (6, -11) {\small $(2, 5)$};

\draw[thick] (8, -6.928) arc (60:0:4)arc  (240:420:4) arc (-60:-120:4); 
\node at (12, -7) {\small $(1,4)$};

\draw[thick] (14, -3.464) arc (-60:120:4) ; 
\node at (12, 0) {\small $(2,4)$};

\draw[thick] (14, 3.464) arc (-60:120:4);  
\node at (12, 7) {\small $(1, 3)$}; 

\draw[thick] (8, 6.928) arc (-60:120:4);  
\node at (6, 11) {\small $(3, 2)$};

\node at (0, 7) {\small (4,2)};
\node at (-6, 3) {\small (2,1)};
\node at (-6, -4) {\small (1, 6)};
\node at (0, -7) {\small (3,5)};
\node at (6, -5) {\small (5, 4)};
\node at (6,5) {\small $(6, 3)$}; 

\draw[thick] (4, 0) arc (-60:120:2) arc (180:-120: 4);
\draw[thick] (4, 0) arc (60:-120:2) arc (-180: 60:4); 
\draw[thick] (-2, 3.464) arc (120:300:4);
\draw[thick] (-2, -3.464) arc (-240:0:4);
\draw[thick] (-4, 0) arc (-60:-120:4) arc (120:300:4)   ;
\draw[thick] (-8, 0) arc (240:60:4);
\draw[thick] (-2, 3.464) arc (240:0:4);

\draw[thick] (2, 3.464) arc (0:180:2);

\node at (0, 0) {$\cdot$};
\node at (0,-1.2) {\bf 0};
\node at (0, -21) {$V_{0, 1, 1}^{*}$};

\end{tikzpicture}
%%%%%%%%%%%%%%%%
%%%%%%%%%%%%%%%%
\quad
%%%%%%%%%%%%%%%%%%
%%%%%%%%%%%%%%%%%%
\begin{tikzpicture}[scale=0.12]
\fill[lightgray] (2, 3.464) arc (60:300:4) arc (-120:60:2) arc (-120:180:4);

\draw[thick] (4, 0) arc (-120:180: 4);
\draw[thick] (4, 0) arc (60:-120:2) arc (-180: 60:4); 
\draw[thick] (-2, 3.464) arc (120:300:4);
\draw[thick] (-2, -3.464) arc (-240:0:4);
\draw[thick] (-4, 0) arc (-60:-120:4) arc (120:300:4)   ;
\draw[thick] (-8, 0) arc (240:60:4);
\draw[thick] (-2, 3.464) arc (240:0:4);

\draw[thick] (2, 3.464) arc (60:300:4);
\draw[thick] (10, -3.464) arc (-120:120:4);
\draw[thick] (8, 6.93) arc (-60:180:4);
\draw[thick] (14, 3.464) arc (-60:120:4);

\node at (0, 7) {\small (2,2)};
\node at (-6, 3) {\small (2,1)};
\node at (-6, -4) {\small (1, 6)};
\node at (0, -7) {\small (3,5)};
\node at (6, -5) {\small (5, 4)};
\node at (12, 0) {\small $(2, 4)$}; 
\node at (12, 8) {\small $(1, 3)$}; 
\node at (6, 10) {\small $(3, 2)$}; 

\draw[thick] (2, 10.392) arc (180:120:4) arc (0:180:4) arc (60:0:4); 
\node at (0, 14) {\small $(1,2)$};

\draw[thick] (-2, 10.392) arc (0: 240:4) ; 
\node at (-6, 11) {\small $(3,1)$}; 

\draw[thick] (-8, 6.928) arc (240:180:4) arc (60:240:4) arc (120:60:4); 
\node at (-12, 7) {\small $(1,1)$};

\draw[thick] (-14, 3.464) arc (120:300:4) ; 
\node at (-12, 0) {\small $(3,6)$};

\draw[thick] (-14, -3.464) arc (120:300:4) arc (180:120:4)   ;
\node at (-12, -7) {\small $(1,6)$};

\draw[thick] (-8, -6.928) arc (120:360:4);  
\node at (-6, -11) {\small $(2, 6)$};

\draw[thick] (-4, -13.856) arc (180:360:4) arc (240:180:4);  
\node at (0, -15) {\small $(1,5)$};

\draw[thick] (2, -10.392) arc (-180: 60:4) ;
\node at (6, -11) {\small $(2, 5)$};

\draw[thick] (8, -6.928) arc (60:0:4)arc  (240:420:4) arc (-60:-120:4); 
\node at (12, -7) {\small $(1,4)$};

\draw[thick] (14, -3.464) arc (-60:120:4) ; 
\node at (12, 0) {\small $(2,4)$};

\node at (0, 0) {\small $\cdot$};
\node at (0,-1.3) { \bf 0};
\node at (0, -20) {$V_{0, 4, 1}^{*}$};

\end{tikzpicture}
%%%%%%%%%%%%%%%%%%%%
%%%%%%%%%%%%%%%%%%%%

\end{center}
\caption{More about configurations around the unit circle : $V_{0, 1,1}^{*}$ and $V_{0, 4, 1}^{*}$}
\end{figure}
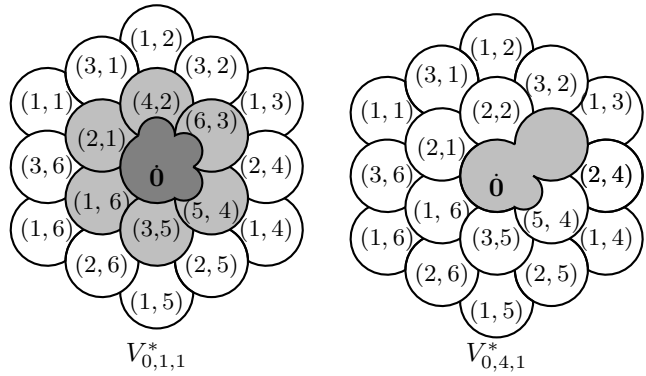

%%%%%%%%%%%%%%%
%%%%%%%%%%%%%%%
\begin{lem}  \label{lem-4-3}
For any $(k, \ell)$, $1 \le k, \ell \le 6$, $V_{k,\ell}^{*} \subset V_{0, k,\ell}^{*}$. 
\end{lem}
%%%%%%%%%%%%%%%
%%%%%%%%%%%%%%%
{\bf Proof.}  From Lemma \ref{lem-4-1}, we see $(z, \infty) \in (z, V_{0, k, \ell}^{*})$ 
whenever $z \in V_{k, \ell}$.  Then Lemma \ref{lem-4-2} shows the claim of 
this lemma. 
\qed\\

%%%%%%%%%%%%%%%
%%%%%%%%%%%%%%%
\begin{cor} \label{cor-4-1}
For any $(k, \ell)$, we have 
\[
V_{0, k, \ell}^{*} = 
\bigcup_{\alpha} \left(V_{0, k', \ell'}^{*}\right)^{-1} - \alpha 
\]
where $\alpha \in \eta \cdot \mathfrak o(\sqrt{-3})$ such that there exists $z \in 
V_{0, k', \ell'}$ with $T(z) \in V_{0, k, \ell}$.  
In particular, any pair of the interior of the right side element are disjoint.  Moreover,  
$V_{0, k, \ell} \cap \{z : |z|\ge \sqrt{3} + 1\} = \{z : |z|\ge \sqrt{3} + 1\}$.  
\end{cor}
%%%%%%%%%%%%%%%
%%%%%%%%%%%%%%%
\begin{prop} \label{prop-4-1}
For any $(k, \ell)$, we have $V_{k, \ell}^{*} = V_{0, k, \ell}^{*}$. 
\end{prop}
{\bf Proof.}   We fix $(k, \ell)$ and start with any $w \in (V_{0, k, \ell}^{*})^{\circ}$.  
From Lemma \ref{lem-4-2}, there exists a unique $\alpha_{1} \in \mathcal{J}$ and 
$(k', \ell')$ such that  $w \in \alpha_{1} + (V_{0, k', \ell'})^{-1}$.  Hence $w' = w + \alpha_{1} 
\in (V_{0, k', \ell'}^{*})^{-1}$.  This means that for $z \in V_{k, \ell}$ we have 
$\hat{T}^{-1}(z, w) = \left(\tfrac{1}{z - \alpha_{1}}, \tfrac{1}{w'}\right)$.  Since 
$\cup_{k, \ell} \partial V_{0, k, \ell}$ is of Lebesgue measure $0$, we can assume 
that $w' \in (V_{0, k', \ell'}^{*})^{\circ}$.  Inductively, we can choose a sequence 
$(\alpha_{1}, \alpha_{2}, \ldots)$ such that 
\[
z^{(n)} =\confrac{1}{\alpha_{n}} + \confrac{1}{\alpha_{n-1}} + \cdots + 
\confrac{1}{\alpha_{2}} + \confrac{1}{\alpha_{1}+z} \in U
\]
for any $n \ge 1$.  Moreover, it is easy to show that 
\[
- \frac{1}{w} = \confrac{1}{\alpha_{1}} + \confrac{1}{\alpha_{2}} + \cdots + 
\confrac{1}{\alpha_{n}} + \cdots .  
\]
Then we see 
\[
 \confrac{1}{\alpha_{1}} + \confrac{1}{\alpha_{2}} + \cdots + 
\confrac{1}{\alpha_{n}}
= \frac{q_{n-1}(z^{(n)})}{q_{n}(z^{(n)})}, 
\]
which shows that $w \in V_{k, \ell}^{*}$.  
\qed \\ 
%%%%%%%%%%%%%%%%%%%%%%
%%%%%%%%%%%%%%%%%%%%%% 

From Proposition \ref{prop-4-1}, we see $\hat{U} = 
\cup_{k=1}^{6}\cup_{\ell = 1}^{6} V_{k, \ell} \times V_{k, \ell}^{*}$.
%%%%%%%%%%%%%%%%%%%%%%%%%%%%%%%%
\begin{thm} \label{thm-4-1}
$(\hat{U}, \hat{T}, \hat{\mu})$ is a representation of the natural extension of $T$, 
where $\hat{\mu}$ is given by $C_{0} 
\tfrac{d(\lambda \times \lambda)}{|z - w|^{4}}$ and $C_{0}$ is the normalizing constant. 
\end{thm}
\begin{rem} The natural extension means that $\hat{T}$ is the (unique) invertible 
minimal extension of $T$.  It is not hard to see that $\int_{\hat{U}} \frac{d(\lambda \times \lambda)}{|z - w|^{4}} < \infty$. The absolutely continuous measure $\mu$ for $T$ is given by the density function 
\[
h(z) =C_{0} \int_{V_{k, \ell}^{*}} \frac{d\lambda(w)}{|z - w|^{4}} 
\]
for $z \in V_{k, \ell}$.  

\end{rem}
{\bf Proof.} It is obvious that $\hat{T}$ is an extension of $T$ since the first 
coordinate 
of $\hat{T}(z, w)$ is $T(z)$.  From Lemma \ref{lem-4-2}, it is invertible (a.e.).  
The proof of Lemma \ref{lem-4-2} implies that $\hat{T}$ separates the second 
coordinate of $(z, w)$ which implies the minimality.  
\qed \\
%%%%%%%%%%%%%%%
%%%%%%%%%%%%%%%
As a corollary, we have the following. 
%%%%%%%%%%%%%%%%%%
%%%%%%%%%%%%%%%%%%
\begin{cor}
For almost every $z \in U$, we have 
\[
\lim_{n \to \infty}\frac{1}{n} \log |q_{n}(z)| 
= \int_{\hat{U}} \log |w| \, d\hat{\mu}(z, w) 
\]
and the right side is finite and positive. 
\end{cor} 
%%%%%%%%%%%%%%%%%%
%%%%%%%%%%%%%%%%%%
{\bf Proof.}  For any $(z, w) \in \hat{U}$, we write $\hat{T}^{n}(z, w) = 
(z_{n}, w_{n})$.  Since $\hat{T}$ is contractive along the $w$-direction, we see 
that $\left|w_{n} - \tfrac{q_{n}(z)}{q_{n-1}(z)}\right| \to 0$.  Hence 
\begin{equation}    \label{eq-4-3}
\lim_{N \to \infty} \frac{1}{N} \sum_{n=1}^{N} \log |w_{n}| 
=  \lim_{N \to \infty} \frac{1}{N} \sum_{n=1}^{N} 
\log |\tfrac{q_{n}(z)}{q_{n-1}(z)}|.   
\end{equation} 
The left side converges for almost every $(z, w)$ by the ergodic theorem and 
the ergodicity of $\hat{T}$.  Hence, for almost every $z$, 
there exists $w_{0}$ such that $(z, w_{0}) \in \hat{U}$ by Fubini's theorem and the limit of 
\eqref{eq-4-3} is 
equal to $\int_{\hat{U}} \log|w|d\mu(w)$, where $\log |w|$ is integrable 
with respect to $\mu$.   Once we have $(z, w_{0})$, \eqref{eq-4-3} shows $(z, \infty)$ 
has also the same limit.  Obviously, the right side is $\lim \tfrac{1}{N} \log |q_{N}(z)|$. 
\qed \\

%%%%%%%%%%%%%%%
%%%%%%%%%%%%%%%
\section{Monotonicity of $|q_{n}|$}
In the previous section, we have shown that $\hat{T}$ is defined on $\hat{U}$ into itself.  It is easy to see that if $(b_{1}(z), \ldots, b_{n}(z))$ is admissible for any 
$n \ge 1$, then $|q_{n+1}(z)| > |q_{n}|$ holds for any $n \ge 1$.   This is basically 
included in the definition of $(\hat{U}, \hat{T})$ as $\hat{T}^{n+1}(z, \infty) = 
\left(T^{n}(z), -\tfrac{q_{n+1}}{q_{n}}\right)$.   
However, if there exists $m_{0}$ such that $(b_{1}(z), \ldots, b_{m_{0}}(z))$ is not 
admissible, then we have to look at some of the boundaries of $V_{k, \ell}$ to 
show $|q_{m+1}(z)| > |q_{m}|$ for $m \ge m_{0}$.  We note that  
$(b_{1}(z), \ldots, b_{m}(z))$ is not admissible for $m \ge m_{0}$.  Also we have to 
check the cases $|T^{m_{0}}(z)|= -\zeta$ and $= \overline{\zeta}$ because 
the expansion of $z$ is not induced from $T$.  

In this section, we give a brief explanation of the monotonicity of $(q_{n}(z) : n \ge 1)$ for any $z \in U$.  We prove the following: 
%%%%%%%%%%%%%%%%
%%%%%%%%%%%%%%%%
\begin{thm} \label{thm-5-1}
For any $z \in U$, if $T^{n}(z) \ne 0$ for $n \ge 1$, then $|q_{n}(z)|<|q_{n+1}(z)|$ 
holds.
\end{thm}
The rest of this section is the proof of Theorem \ref{thm-5-1} for $z \in U$ such that 
there exists $n \ge 1$ such that $(b_{1}(z), \ldots, b_{n}(z))$ is not 
admissible.  

 If $T^{n}(z)=0$, then it is obvious that $z = \tfrac{p_{n}(z)}{q_{n}(z)}$ and 
we can ignore $\tfrac{p_{m}}{q_{m}}$ for $m > n$.  In the following, we always 
assume that $T^{n}(z) \ne 0$ and then $b_{n+1}(z) \ne \infty$.   Note that 
$T^{n}(z) \ne 0$ implies $z \ne 0$. 

We see from Definition 1 and Proposition \ref{prop-4-1} the following. 
%%%%%%%%%%%%%%%%
%%%%%%%%%%%%%%%%
\begin{lem} \label{lem-5-1}
We have \\
{\rm (i)}  $V_{2, \ell}^{*} = V_{3, \ell -1}^{*}$, where we read $1 -1 = 6$ (mod $6$), 
\\
{\rm (ii)} $\left( \bigcap_{k=1}^{6} \{z \in \mathbb C : |z - \eta_{k}|> 1\} \cap
\{z: |z|> 1\}\right)\subset 
V_{k, \ell}^{*}$\\
for any $1 \le k, \ell \le 6$.  
\end{lem} 
{\bf Proof. } From Definition 1, it is easy to see that 
$V_{0, 2, \ell}^{*} = V_{0, 3, \ell -1}^{*}$ and Proposition \ref{prop-4-1} shows the 
assertion of (i), see Fig. 4.   The assertion (ii) also follows the same 
reasoning, see Fig. 7.   \qed \\

There are two cases \\
(a) $T^{n}(z) \notin \{-\zeta, \overline{\zeta}\}$ for any $n \ge 1$, \\
(b) $T^{n}(z) \in \{-\zeta, \overline{\zeta}\}$ for some $n \ge 1$ and $|T^{n-1}(z)| 
< 1$.

First, we discuss the case (a).   We consider 
$T^{n}(\langle a_{1}, a_{2}, \ldots , a_{n} \rangle) \cap U^{\circ}$, but for the 
simplicity, we omit writing ``$\cap U^{\circ}$".  

A sequence $(a_{1}, a_{2}, \ldots, a_{n})$ in $\mathcal J^{n}$, $n\ge 1$, 
is said to be non-admissible if $\langle  a_{1}, a_{2}, \ldots, a_{n} \rangle$ is not 
empty but its interior is empty.  For any $\alpha \in \mathcal J \setminus 
\{0\}$, $\langle \alpha \rangle$ has an inner point.  Thus,  for $z \in U$, 
if there exists $n \ge 2$ such that $(b_{1}(z), b_{2}(z), \ldots , b_{n}(z))$ is 
non-admissible, then there exists $n_{*} \ge 2$ such that 
$(b_{1}(z), b_{2}(z), \ldots , b_{n_{*}}(z) )$ is admissible and 
$( b_{1}(z), b_{2}(z), \ldots , b_{n_{*}}(z), b_{n_{*}+1}(z) )$ is non-admissible.   
\begin{prop} \label{prop-5-2}
For $z \in U$ such that $n_{*}$ exists, 
$T^{n_{*}+1}(\langle b_{1}(z), b_{2}(z), \ldots , b_{n_{*}}(z), b_{n_{*}+1}(z) \rangle)$ 
is one of the following $L_{j}$, $1 \le j \le 3$, see Fig.8: 
\begin{figure} \label{fig-5-1}
\begin{center}
\begin{tikzpicture}[scale=0.75]
\draw[thick] (2, 3.464)--(-2, 3.464); 

\draw[thick] (2, -3.464)--(4, 0); 

\draw[thick] (-4, 0)--(-2, -3.464);

\draw[thick] (2, 3.464) arc (-30:-150:2.31);
\draw[thick] (-4, 0) arc (90:-30:2.31);
\draw[thick] (2, -3.464) arc (210:90:2.31);

\draw[thick] (4, 0) -- (-4, 0);
\draw[thick] (2, 3.464) -- (-2, -3.464); 
\draw[thick] (-2, 3.464) -- (2, -3.464); 

\draw[thick] (4, 0) arc (30:270: 2.31);
\draw[thick](2, 3.464) arc (30:-210:2.31);
\draw[thick] (-4, 0) arc (150:-90:2.31);

\node at (0, 0) {\Large $\cdot$}; 
\node at (-0.4, 0) {\Large $O$}; 
\node at (4.2, -0.3) {$1$}; 
\node at (3.4, 3.52) {$\zeta = \frac{1 + \sqrt{-3}}{2}$} ; 
\node at (0, 4) {$L_{1}$}; \node at (-3.5, -2) {$L_{2}$}; \node at (3.5, -2) {$L_{3}$};
\node at (1.5, 2) {$L_{4}$}; \node at (0, 2) {$L_{7}$}; \node at (-2, 2) {$L_{10}$};

\end{tikzpicture}
\end{center}
\caption{non-admissible cylinder sets}

\end{figure}
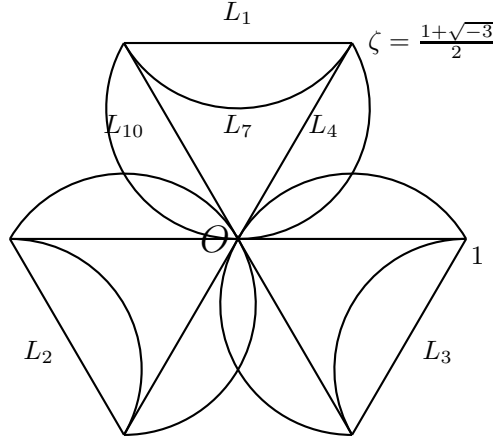
%%%%%%%%%%%%%%%%%%%%%
%%%%%%%%%%%%%%%%%%%%%
%%%%%%%%%%%%%%%%%%%%%
\[
\begin{array}{l}
L_{1}=\{z = t \cdot \zeta + (1-t)\cdot (-\overline{\zeta}) : 0 < t < 1 \} \\ 
L_{2}=\{z =  t \cdot (-\zeta) + (1-t) \cdot(- 1) : 0 < t < 1 \}\\
L_{3} = \{z = t \cdot 1 + (1-t)\cdot \overline{\zeta} : 0 < t < 1 \}  
\end{array} . 
\]
\end{prop}
{\bf Sketch of the Proof. }
We can check this by cheking all possible choices of $n_{*}$ with  
$T^{n_{*}}(\langle a_{1}, a_{2}, \ldots a_{n} \rangle ) = U_{k, \ell}$.
\qed\\

Now we can give all possible choices of $T^{n}(\langle a_{1}, a_{2}, \ldots a_{n} 
\rangle )$ for non-admissible sequences.   Then it is not hard to show that the 
monotonicity of $|q_{n}|$ for non-admissible sequences.
%%%%%%%%%%%%%%%%
%%%%%%%%%%%%%%%%
\begin{prop} \label{prop-5-3}
For any non-admissible sequence $(a_{1}, a_{2}, \ldots, a_{n})$, 
$T^{n}(\langle a_{1}, a_{2}, \ldots, a_{n} \rangle )$ is one of $L_{j}$, $1 \le j \le 12$, 
where $L_{j}$, $4 \le j \le 12$ are given as 
\[
\begin{array}{lll}
 L_{4}&=&\{z = t \zeta + (1-t)(- {\zeta}) : 0 < t < 1 \} \\ 
L_{5} &= &\{z =  t\cdot (-1)+(1-t) \cdot 1 : 0 < t < 1 \}\\ 
L_{6}&=&\{z = t \overline{\zeta} + (1-t)(-\overline{\zeta}) : 0 < t < 1 \} \\
L_{7}& =& \{z \in U^{\circ} : \left| z - \tfrac{2}{3} \eta_{2} \right| < \tfrac{\sqrt{3}}{3} \} \\
L_{8} &= &\{z \in U^{\circ} : \left| z - \tfrac{2}{3} \eta_{4} \right| < \tfrac{\sqrt{3}}{3} \} \\
L_{9}&=&\{z \in U^{\circ} : \left| z - \tfrac{2}{3} \eta_{6} \right| < \tfrac{\sqrt{3}}{3} \\
L_{10}& =& \{z \in U^{\circ} : \left| z - \tfrac{1}{3} \eta_{2} \right| < \tfrac{\sqrt{3}}{3}\} \\
L_{11}&=&\{z \in U^{\circ} : \left| z - \tfrac{1}{3} \eta_{4} \right| < \tfrac{\sqrt{3}}{3} \} \\
L_{12}& = & \{z \in U^{\circ} : \left| z - \tfrac{1}{3} \eta_{6} \right| < \tfrac{\sqrt{3}}{3}\}
\end{array} . 
\]
\end{prop}
{\bf Proof.}  We already saw that $L_{j}$, $j=1, 2, 3$ are possible choices, see 
Proposition \ref{prop-5-2}. 
Now suppose that $T^{m}(\langle a_{1}, a_{2}, \ldots, a_{m} \rangle) = L_{1}, 
L_{2}$ or $L_{3}$ for $(a_{1}, a_{2}, \ldots, a_{m})$.  Then each $(m+1)$th 
coordinate is uniquely determined as $\eta_{5}, \eta_{3}$ $ \eta_{1}$ and 
\begin{equation} \label{eq-5-1}
\begin{array}{l}
T^{m+1}(\langle a_{1}, a_{2}, \ldots, a_{m},  \eta_{5} \rangle) = L_{7}, \\
T^{m+1}(\langle a_{1}, a_{2}, \ldots, a_{m},  \eta_{3} \rangle) = L_{9},\\
T^{m+1}(\langle a_{1}, a_{2}, \ldots, a_{m}, \eta_{1} \rangle) = L_{8},  
\end{array}
\end{equation}
respectively.  Moreover, 
we have 
\begin{equation} \label{eq-5-2}
\begin{array}{l}
T^{m+2}(\langle a_{1}, a_{2}, \ldots, a_{m},  \eta_{5}, \eta_{5} \rangle) = L_{10}, \\
T^{m+2}(\langle a_{1}, a_{2}, \ldots, a_{m},  \eta_{3}, \eta_{1} \rangle) = L_{11}, \\
T^{m+2}(\langle a_{1}, a_{2}, \ldots, a_{m}, \eta_{1}, \eta_{3} \rangle) = L_{12}.
\end{array}
\end{equation}
Again, each $(m+2)$th coordinate is uniquely determined. 
Once we have one of these $12$, it is easy to see that no other line segment or 
arc appears as an image of a cylinder set arising from a non-admissible sequence. 
\qed\\
%%%%%%%%%%%%
%%%%%%%%%%%%
\begin{prop} \label{prop-5-4} 
For any $z \in U$, if $T^{n}(\langle b_{1}(z), b_{2}(z), \ldots, b_{n}(z) \rangle)$ 
is non-admissible, then $|q_{n-1}(z)| <|q_{n}(z)|$. 
\end{prop}
%%%%%%%%%%%%
%%%%%%%%%%%%
{\bf Proof.}  
There are two cases: (i) $n = n_{*} + 1$, (ii) $n > n_{*}+1$.  First, we prove the case (i). 
As we discussed in the proof of Proposition \ref{prop-5-2}, there are two cases.  

In the first case, we always have $|b_{n_{*}+1}(z)| \ge 2 \sqrt{3}$ and so 
$-\tfrac{q_{n_{*}+1}(z)}{q_{n_{*}}(z)}\in \left(U_{k, \ell}\right)^{\circ}$ for any 
$(k, \ell)$, see Lemma \ref{lem-5-1}(ii).  In particular, if 
\[
T^{n_{*}+1}(\langle b_{1}(z), b_{2}(z), \ldots, b_{n_{*}}(z), b_{n_{*}+1}(z) \rangle) 
= L_{1}, L_{2} \,\, \mbox{or}\,\, L_{3}, 
\]
then  $-\tfrac{q_{n_{*}+1}(z)}{q_{n_{*}}(z)}\in \left(V_{6, 2}^{*}\right)^{\circ}, 
\left(V_{6, 4}^{*}\right)^{\circ}$ or $\left(V_{6, 6}^{*}\right)^{\circ}$ respectively.

Now, we consider the second case.  Again, we follow the proof of Proposition 
\ref{prop-5-2}.  If 
\[
T^{n_{*}}(\langle b_{1}(z), b_{2}(z), \ldots, b_{n_{*}}(z) \rangle) 
= U_{3, 1}, U_{4, 1} \,\, \mbox{or}\,\, U_{5, 4}, 
\]
then we see $-\tfrac{q_{n_{*}}(z)}{q_{n_{*}-1}(z)} \in (V_{3,4}^{*})^{\circ}$ or 
$(V_{2, 2}^{*})^{\circ}$.  
From Lemma \ref{lem-5-1}(i), we have 
\[
-\frac{q_{n_{*}}(z)}{q_{n_{*}-1}(z)} \in (V_{2, 5}^{*})^{\circ} \,\,\mbox{or} \,\, 
(V_{3, 1}^{*})^{\circ}  . 
\]
If $|b_{n_{*} + 1}(z)| > \sqrt{3}$, then it is obvious that we have 
$-\tfrac{q_{n_{*}+1}(z)}{q_{n_{*}}(z)}\in (V_{6, 4}^{*})^{\circ}$, the same as before.  
If $b_{n_{*} + 1}(z) = \eta_{2}$ or $\eta_{6}$,  then we can find $z' \in V_{2, 5}$ or 
$V_{3, 1}$ such that $\left[\tfrac{1}{z'}\right] = \eta_{2}$ or $\eta_{6}$, respectively. 
and then $T(z') \in V_{6,4}$.  Hence, we get the same conclusion: 
\[
  -\frac{q_{n_{*}+1}(z)}{q_{n_{*}}(z)}\in (V_{6, 4}^{*})^{\circ} .  
\]
By the same way, we have the following: if 
\[
T^{n_{*}}(\langle b_{1}(z), b_{2}(z), \ldots, b_{n_{*}}(z) \rangle)  =  \left.
\begin{array} {l} 
U_{3, 3}, U_{4, 3}, U_{5, 6}\\
U_{3, 5}, U_{4, 5}, U_{5, 2}
\end{array}\right\}, \,\,  \mbox{then} \,\,  -\frac{q_{n_{*}+1}(z)}{q_{n_{*}}(z)}
\in  \left\{
\begin{array}{l} 
(V_{6, 2}^{*})^{\circ}\\ (V_{6, 6}^{*})^{\circ} \end{array}\right. , 
\]
respectively.

Now we consider the case (ii).  From the proof of Proposition \ref{prop-5-3}, 
we have cycles 
\[
L_{1}\to L_{4} \to L_{7} \to \left\{ \begin{array}{ll}  L_{1}\\ L_{10} \to L_{1} \\
                                                                                                L_{10} \to L_{10} \cdots
\end{array}\right.
\]
and 
\[
-\tfrac{q_{n+1}(z)}{q_{n}(z)} \in  
\left\{\begin{array}{lll}
 (V_{6,2}^{*})^{\circ} & \mbox{for} & L_{1}, \\
 (V_{4,2}^{*})^{\circ} or (V_{5,2}^{*})^{\circ} & \mbox{for} & L_{4}, \\
 (V_{6, 3}^{*})^{\circ} or (V_{3, 3}^{*})^{\circ} & \mbox{for} & L_{7}, \\  
 (V_{2, 1}^{*})^{\circ} or (V_{3, 3}^{*})^{\circ} & \mbox{for} & L_{10}
\end{array}\right.
\]
by the 
same reasoning.  The other lines and arcs are the same.  Thus the basic lemma 
shows the conclusion.    
\qed\\

Finally, we consider the case (b). 
If $T^{n}(z) = -\zeta$ or $\overline{\zeta}$ for some $n \ge 0$, then we have to define the expansion of $z$ according to the definition of the expansions of 
$-\zeta$ and $\overline{\zeta}$.  If $n=0$, then 
it is already shown in Lemma \ref{lem-2-1}.  Thus, we consider the case $n \ge 1$ 
with $T^{m}(z) \ne -\zeta, \overline{\zeta}$ for $0 \le m < n$.  
We recall the expansions of $-\zeta$ and $\overline{\zeta}$ in 
Lemma \ref{lem-2-1}. 
Concerning \eqref{eq-2-6} and \eqref{eq-2-7}, we can choose sequences 
$(z_{-\zeta, m} : m \ge 1)$, $(z_{\overline{\zeta}, m} : m \ge 1)$ and 
$(m_{n}: n \ge 1)$ such that  $z_{-\zeta, m}, z_{\overline{\zeta}, m} \in V_{6,5}$ 
and 
\[ 
\left\{ \begin{array}{lll}
b_{n}(z_{-\zeta, m}) &=& b_{n}(-\zeta), \\
b_{n}(z_{\overline{\zeta}, m})& = &b_{n}(\overline{\zeta}) 
\end{array} \right.
\]
for $m \ge m_{n}$, respectively.  Apparently, these sequences converge to $-\zeta$ 
and $\overline{\zeta}$, respectively.  

%%%%%%%%%%%%%%%%%%%%%%
\begin{lem} \label{lem-5-3}
If $T^{n}(z) = - \zeta$ or $\overline{\zeta}$ and $|T^{n-1}(z)|<1$, then 
$-\tfrac{q_{n}(z)}{q_{n-1}(z)} \in \left(V_{6,5}^{*}\right)^{\circ}$.  
\end{lem}
%%%%%%%%%%%%%%%%%%%%%%
{\bf Proof.} It would be enough to show that 
the first case, i.e., $T^{n}(z)= -\zeta$ and $|T^{n-1}(z)|<1$. 
If $|b_{n}(z)|> \sqrt{3}$, then, from Lemma \ref{lem-5-1}(i), we see 
$-\tfrac{q_{n}(z)}{q_{n-1}(z)} \in \left(V_{6, 5}^{*}\right)^{\circ}$.  
So we consider the case $|b_{n}(z)|= \sqrt{3}$.   
There are four possibilities: \\
\begin{equation}
\frac{1}{T^{n-1}(z)} = \left\{ \begin{array}{l}
-2 \\ -2\zeta -1 \\ -2\zeta + \overline{\zeta} \\ 2\overline{\zeta} 
\end{array} \right.   .  
\end{equation}
These imply 
\begin{equation}
b_{n}(z) = \left\{ \begin{array}{l}
\eta_{3} \\\eta_{4} \\ \eta_{5} \\ \eta_{6} \end{array}\right.  ,  
\end{equation}
respectively. 
First, we consider the admissible sequence $(b_{1}(z), b_{2}(z), \ldots, b_{n-1}(z))$.  
In this case, one can easily show the conclusion, since 
these four imply $T^{n-1}(z)$ ``$\in V_{2, 4}$ or $V_{3, 3}$", ``$\in V_{2, 3}$", 
``$\in V_{3,2}$" or ``$\in V_{3, 1}$ or $V_{2, 2}$", respectively.
These show $-\tfrac{q_{n}(z)}{q_{n-1}(z)} \in \left(V_{6, 5}^{*}\right)^{\circ}$ 
for all cases.  
Next we consider the non-admissible $(b_{1}(z), b_{2}(z), \ldots, b_{n-1}(z))$.  
It is easy to see that this occurs only when $\frac{1}{T^{n}(z)} = -2$ or 
$2 \overline{\zeta}$.  For these cases, we can reduce it to the admissible case 
as we did in the proof of Proposition \ref{prop-5-4}.  
Once we get $T^{n}(z) = - \zeta$ with $-\tfrac{q_{n}(z)}{q_{n-1}(z)} \in 
\left(V_{6,5}^{*}\right)$, we have 
\[
\begin{array}{lll}
\frac{1}{-\zeta} - \eta_{2} = -\zeta & \mbox{with}& \frac{q_{n}(z)}{q_{n-1}(z)} \in 
\left(V_{2, 5}\right)^{\circ} \\
\frac{1}{-\zeta} - \eta_{2} = -\zeta & \mbox{with}& \frac{q_{n}(z)}{q_{n-1}(z)} \in 
\left(V_{6, 4}\right)^{\circ} \\
\frac{1}{-\zeta} - \eta_{3} = 1 & \mbox{with}& \frac{q_{n}(z)}{q_{n-1}(z)} \in 
\left(V_{3, 6}\right)^{\circ} \\
\mbox{and} & & \\
\frac{1}{1} - \eta_{1} = -\zeta & \mbox{with}& \frac{q_{n}(z)}{q_{n-1}(z)} \in 
\left(V_{6, 5}\right)^{\circ} 
\end{array}  .  
\]
Then it goes on cyclically, see (8) and Fig. 6.

Thus we get the assertion of this lemma. \qed\\

Finally, following the idea of the proof of Proposition \ref{prop-5-4}, we have 
%%%%%%%%%%%
\begin{prop}  \label{prop-5-5}
If $T^{n}(z) = - \zeta$ or $\overline{\zeta}$ and $|T^{n-1}(z)|<1$, then 
$|q_{m}(z)|< |q_{m+1}(z)|$ for any $m \ge n$. 
\end{prop}
%%%%%%%%%%%
From all of the above discussion, we conclude that Theorem \ref{thm-5-1} 
holds. 
\\[5pt]
 
{\bf Acknowledgments.}\\ 
The authors would like to thank Hiromi Ei for her helpful suggestions. 
This research was partially supported by JSPS grants 
24K06785 (the first author), and 20K03642 (the second author).

%%%%%%%%%%%%%%%
%%%%%%%%%%%%%%%

\end{document}